\documentclass[twocolumn]{autart}    

\usepackage{graphicx}          
                               
\usepackage{amsmath}    		
\usepackage{graphics} 		
\usepackage{setspace}		
\usepackage{longtable}          
\usepackage{verbatim}
\usepackage{fancyhdr}      
\usepackage{bm}
\usepackage[noadjust]{cite}

\usepackage{color}

\usepackage{pdfpages}
\usepackage{appendix}

\usepackage{epsfig} 
\usepackage{mathptmx} 
\usepackage{times} 
\usepackage{amsmath} 
\usepackage{amssymb}  
\usepackage{amsfonts}
\usepackage{tikz,tabstackengine,amsmath}
\usetikzlibrary{shapes,shadows,arrows}
\usetikzlibrary{shapes,shapes.geometric,arrows.meta,fit,calc,positioning,automata}
\usepackage{pgfplots}
\usepgfplotslibrary{groupplots}
\usepackage{subfig}
\usepackage{mathtools}
\usepackage{adjustbox}
\usepackage{xcolor}
\usepackage{color}
\usepackage{tabularx}
\usepackage [english]{babel}
\usepackage [autostyle, english = american]{csquotes}
\usepackage{caption}
\usepackage{todonotes}
\usepackage{dsfont}
\usepackage[utf8]{inputenc}

\usepackage{enumitem}
\newtheorem{remark}{Remark}
\newtheorem{obs}{Observation}
\newcommand{\mb}{\mathbb}
\newcommand{\mc}{\mathcal}

\newcommand{\eps}{\epsilon}

\newcommand{\hf}{\tfrac{1}{2}}
\newcommand{\mf}{\mathfrak}
\newcommand{\trns}{\intercal}
\newcommand{\dsR}{{\mathds{R}}}

\newcommand{\tex}{{\text{exp}}}

\newcommand{\mfN}{{\mathfrak{N}}}
\newcommand{\mfK}{{\mathfrak{K} }}
\newcommand{\mfT}{{\mathfrak{T} }}

\definecolor{myGrn}{rgb}{0,0.5,0}

\usepackage[normalem]{ulem}

\begin{document}

\begin{frontmatter}

\title{Exploratory LQG Mean Field Games with Entropy Regularization \thanksref{footnoteinfo}} 

\thanks[footnoteinfo]{This work was presented in the CMS Winter Meeting, December 6-9, 2019, Toronto, Canada, in the 11th Workshop on Dynamic Games in Management Science, October 24-25, 2019, Montreal, Canada, and in the virtual SIAM conference on Financial Mathematics \& Engineering, June 1-4, 2021. SJ and DF would like to, respectively, acknowledge support from the Natural Sciences and Engineering Research Council of Canada (grants RGPIN-2018-05705 and RGPAS-2018-522715) and the Fonds de Recherche du Qu\'ebec--Nature et Technologies (funding reference number 258061). Corresponding author Dena Firoozi.}

\author[Paestum]{Dena Firoozi}\ead{dena.firoozi@hec.ca},    
\author[Rome]{Sebastian Jaimungal}\ead{sebastian.jaimungal@utoronto.ca}               

\address[Paestum]{Department
of Decision Sciences, HEC Montr\'eal, Montreal,
QC, Canada}  
\address[Rome]{Department of Statistical Sciences, University of Toronto, ON, Canada}             
          
\begin{keyword}                           
LQG mean field games; entropy-regularization; exploration.              
\end{keyword}                             

\begin{abstract}  
We study a general class of entropy-regularized multi-variate LQG mean field games (MFGs) in continuous time with $K$ distinct sub-population of agents. We extend the notion of control actions to control distributions (exploratory actions), and explicitly derive the optimal control distributions for individual agents in the limiting MFG. We demonstrate that the optimal set of control distributions yields an $\epsilon$-Nash equilibrium for the finite-population entropy-regularized MFG. Furthermore, we compare the resulting solutions with those of classical LQG MFGs and establish the equivalence of their existence. 
\end{abstract}

\end{frontmatter}

\section{Introduction}
Mean Field Game (MFG) systems theory establishes the existence of approximate Nash equilibria and the corresponding individual strategies for stochastic dynamical systems in games involving a large number of agents. The equilibria are termed $\epsilon$-Nash equilibria and are generated by the local, limited information control actions of each agent in the population. The control actions constitute the best response of each agent with respect to the behaviour of the mass of agents. Moreover, the approximation error, induced by using the MFG solution, converges to zero as the population size tends to infinity. The analysis of this set of problems originated in \cite{HuangCDC2003b, HuangCIS2006, HuangTAC2007}, and independently in \cite{Lasry2006a, Lasry2006b, Lasry2007} using forward-backward PDEs, known as the analytic approach in the literature. Several approaches have been developed in the past years, principally the probabilistic approach \cite{carmona2018probabilistic}, and the master equation approach \cite{MasterEq2019}. Moreover, many extensions and generalizations exist including but not limited to MFGs with a major agent \cite{Huang2010, NourianSiam2013,CarmonaWang2017,LasryLions2018,CardaliaguetCirant2018,BensoussanSICON2017}, partial observations \cite{SenSiam2016,FirooziCainesTAC2020,BensousanPO2020}, stopping and switching times \cite{FPC-arXiv2018,Tankov2020}, and mean-field type problems \cite{ BensoussanBook2013}. 

MFGs have found numerous applications in engineering problems such as cellular network optimization \cite{aziz2016mean} and coordination of loads in smart grids \cite{KIZILKALE2019} (see \cite{Tembine2017} for a set of interesting applications), and in particular in mathematical finance and economics for characterizing equilibrium price and market equilibria (see \cite{shrivats2020mean,FirooziISDG2017,Gomes2016, CarmonaRev2020} and the references therein) -- to name a few. 

Entropy regularization has been used to develop iterative schemes for solving PDEs. \cite{FPK-entropy} constructs a discrete-time iterative scheme whose solutions converge to that of a class of Fokker-Planck equations. This scheme, for which the time step is governed by the Wasserstein metric on probability measures, views the Fokker-Planck equation as a steepest descent of the negative entropy. \cite{Gomes2007} uses entropy regularization and a lifting to the space of control distributions to develop a discrete-time iterative scheme that approximates the viscosity solutions to the Hamilton-Jacobi-Bellman equation. In the same work, the connection between a similar iterative scheme and an entropy regularized version of the stationary Mather problem is established. It is shown that the unique minimizing measure can be explicitly given in terms of the fixed-point solution to the scheme. Moreover, this solution is related to the difference between the solutions to a backward and a forward fixed-point schemes. In a recent work \cite{XYZ-Langevin2020} entropy regularization is used for designing an endogenous temperature process for Langevin algorithm in non-convex optimization problems.  

Recently entropy regularized stochastic control problems \cite{Neu2017, Geist2019,Lefebvre2020,XYZRL2020} and stochastic games \cite{Savas2019}, where an entropy term is included in the objective function and the space of actions is the space of measures, have received attention in the literature. This research direction aims to provide a theoretical basis for the policy search techniques used in the field of reinforcement learning (RL), where an agent learns about an unknown environment by trial and error using feedback from their own actions and experiences. The agent's aim is to find a suitable action model that maximizes their total cumulative reward.  When the agent takes random actions it is called exploration and, as it is usually costly, an agent faces the dilemma of exploring new states versus maximizing their reward given their current knowledge. This is referred to as the exploration versus exploitation trade-off. A recent trend in RL includes exploration in the objective function through entropy regularization in discrete time \cite{Ziebart2008, Nachum2018, Fox2016}. 

While most of the literature studies entropy regularization using Markov decision processes, a few studies on the continuous-time stochastic diffusion processes have emerged recently, in particular \cite{XYZ-MV-2020, XYZRL2020, Siska2020}.        



In this paper we study a general class of LQG MFGs which are regularized by Shannon entropy in continuous time. We also consider $K$ distinct sub-populations of agents, where agents within a sub-population have the same dynamical and cost functional parameters. Each agent has a multivariate state and control action processes and only observes their own state at each time instant. Following a similar approach as in \cite{XYZRL2020, XYZ-MV-2020} for single-agent systems, we extend the notion of actions to control distributions, and explicitly derive the optimal control distributions for individual agents in the limiting MFG. We demonstrate that the set of obtained distributions yields an $\epsilon$-Nash equilibrium for the finite-population entropy-regularized MFG. We compare the obtained solutions with those of classical LQG MFGs and establish that the entropy regularized solution exists if and only if the classical version exists.

The optimal control distributions may be interpreted as the optimal policy for an RL agent that achieves the best balance of exploration versus exploitation. Interestingly, for the class of MFGs under consideration the optimal distributions are shown to be Gaussian with mean equal to the optimal control action for classical LQG MFGs. 
The variance of the distribution may be interpreted as a measure of the amount of exploration the agent should engage in, while the mean may be interpreted as the agent exploiting their knowledge.

In practice models may be nonlinear as in \cite{shrivats2020mean}, where optimal behaviour and equilibrium pricing in renewable energy markets have been studied. LQG models, however, are of significant importance from the analytical point of view as  their  solutions may be explicitly characterised. This provides intuition and helps understanding the properties and the nature of the system. Moreover, from the application point of view, LQG models may be used to approximate nonlinear, and in some cases intractable, models. Due to these reasons, LQG models have received significant attention in the literature. Hence, we address LQG MFG models in this manuscript as a first step in  exploratory MFG research.

Two works in the context of MFGs related to the current paper are (i) \cite{SaldiKariksiz2020}, where a class of discrete-time entropy-regularized mean field games are studied, while the explicit control distributions are not derived, a Q-learning algorithm is established to compute approximate equilibria ; (ii) \cite{Zariphopoulou2020}, where the solutions to a class of continuous-time LQG MFGs with entropy regularization in the infinite population limit is studied. The authors consider a scalar mean-reverting dynamics, where the diffusion coefficient is impacted by the control action together with a finite-horizon tracking objective function. Agents are assumed to be homogeneous with the same dynamical and cost functional parameters, and to observe their own state, the population state distribution and the population control distribution at each time instant. In addition to Shannon entropy, they also investigate a combination of Shannon entropy and the cross-entropy. They perform all the analysis in the infinite-population limit and obtain Nash equilibria. Finally they propose a policy-gradient algorithm to estimate the unknown parameters. 

Our work distinguishes itself from \cite{Zariphopoulou2020} in several ways: we (i) consider $K$ sub-populations of agents with a general multi-variate linear dynamics (although the diffusion coefficient is not impacted by the control action) and tracking cost functionals; (ii) consider a local information pattern for each agent which only includes its private state and the initial population state distribution; (iii) derive the mean field consistency equations as a set of differential and algebraic equations generating an infinite-population Nash equilibrium and subsequently establish the $\epsilon$-Nash property for the finite-population case; and (iv) establish that the existence of the solutions to the classical and exploratory LQG MFG systems are equivalent.  



Another related line of research including \cite{CarmonaLauriere2019,XinGuo2019,XuWei2020} establishes a dynamic programming principle for mean-field control problems (without regularization), where the problem is formulated as an MDP on the space of measures, and a Q-learning algorithm is designed for learning the optimal solution.

Although it is not the focus of the current paper, it is worth mentioning that there is a rapidly growing literature on RL in mean field games. Among them are
Q-learning based algorithm for MFG control of coupled oscillators \cite{Meyn2014}, model-free Q-learning and actor critic algorithms for mean-field games \cite{Yang2018,Mguni2018}, reinforcement learning in stationary mean field games \cite{MahajanSubra2019}, deep inverse reinforcement learning for discrete MFGs \cite{Yang2017},  MFG Q-learning algorithm with Boltzmann policy \cite{XinGuoLearningMFG2019},  and policy gradient methods for the continuous-time linear-quadratic mean-field controls and games \cite{Wang2020}.

Throughout this paper we use the terms entropy-regularized MFGs and exploratory MFGs interchangeably. 

The organization of the remainder of the paper is as follows. We first briefly review the results of classical LQG MFG in Section \ref{sec:classicLQGMFG}. Then we introduce entropy-regularized LQG MFGs and discuss their solutions in Section \ref{sec:reg_LQGMFG}. The solutions of classical and entropy reqularized LQG MFGs are compared in Section \ref{sec:equivalence} where we also establish the equivalence of their existence. We then briefly discuss potential applications in learning MFGs in Section \ref{sec:apps}. Finally, we conclude the paper in Section \ref{sec:conclusion} with some remarks on the potential applications of the developed framework.

\section{Overview: Classical LQG MFGs}\label{sec:classicLQGMFG}
\subsection{Finite Populations} 
We consider a system consisting of a large number $N$ of minor agents in $K$ types (or equivalently subpopulations) governed by the dynamics  
\begin{align} \label{MinorDyn}
dx_t^i = (A_k x_t^i + F_k x^{(N)}_t+ H_k u^{(N)}_t+B_k u_t^i + b_k(t))dt + D_k dw_t^i,
\end{align}
for $t \in \mf{T}:= [0,\infty)$, $i \in \mfN := \{1, 2, ..., N\}$. The subscript $k \in \mfK :=\{1,\dots,K\}$ with $\, K\leq N$, denotes the type of a minor agent. Here $x^i_t \in \dsR^n$ is the state, $u^i_t \in \dsR^m$ is the control input, $w =\lbrace w^i_t,\, t\geq 0,\, i \in \mfN \rbrace$ denotes $N$ independent standard Wiener processes in $\dsR^r$, where $w^i$ is progressively measurable with respect to the filtration $\mc{F}^w \coloneqq (\mc{F}_t^{w})_{t\in \mfT}$. Moreover, $x^{(N)}_t := \frac{1}{N} \sum_{i \in \mfN} x^i_t$ and $u^{(N)}_t := \frac{1}{N} \sum_{i \in \mfN} u^i_t$ denote, respectively, the average state and control of the minor agents. All matrices ($A_k$, $F_k$, $H_k$, $B_k$, and $D_k$) in \eqref{MinorDyn} are constant and of appropriate dimension; $b_k(t): \mf{T} \rightarrow \dsR^n$ are deterministic functions of time.
\begin{assum} \label{IntialStateAss}
The initial states $\lbrace x^i_0,~ i \in \mfN \rbrace$ defined on $(\Omega, \mathcal{F}, P)$ are normally distributed, mutually independent and also independent of $\mathcal{F}^{w}$, with $\mathbb{E}[x^i_0]=\xi$ and $\sup_{i} \mathbb{E}[(x^i_0)^\trns x^i_0] \leq c < \infty $ ($c$ independent of $N$). 
\end{assum} 
\subsubsection*{Agents types}
Minor agents are given in $K\le N$ distinct types (or equivalently subpopulations), such that the agents in the same sub-population share the same dynamical and cost functional parameters.   
We define the index set $\mc{I}_k$ as
\begin{equation*}
\mc{I}_k = \lbrace i : \theta_i = \theta^{(k)},~ i \in \mfN \rbrace ,
\end{equation*}
for $k\in\mfK :=\{1,\dots,K\}$ and $\theta^{(k)}\in\Theta$, where $\Theta$ is the parameter set.
The cardinality of $\mc{I}_k$ is denoted by $N_k = |\mc{I}_k|$. Then, $\pi^{N} = (\pi_{1}^{N},...,\pi_{K}^N),~ \pi_k^N = \tfrac{N_k}{N}, ~ k \in \mfK$, denotes the empirical distribution of the parameters $(\theta_1,...,\theta_N)$ obtained by sampling independently of the initial conditions and 
the Wiener processes of all agents.
\begin{assum} \label{ass:EmpiricalDistLimit}
There exists $\pi$ such that $\displaystyle\lim_{N \rightarrow \infty} \pi^N = \pi $ a.s.
\end{assum}
%
Each agent aims to  minimize the cost functional
\begin{multline}\label{minorCost}
J^{N}_i(u^i, u^{-i})
= \mb{E} \bigg[ \int_0^{\infty} e^{-\rho t}\bigg( \hf \Vert x_t^i-y_t\Vert^2_{Q_k} + \hf \Vert u_t^i\Vert^2_{R_k}\\+ (x_t^i-y_t)^\intercal S_k u_t^i + \eta_k^\trns (x_t^i-y_t) + n_k^\trns u_t^{i} \bigg) dt \bigg],
\end{multline}
where $y_t = \psi_k x^{(N)}_t$, $u^{-i} = (u^{1}, \dots, u^{i-1}, u^{i+1}, \dots u^{N})$, $\Vert x_t^i-y_t\Vert^2_{Q_k}=(x_t^i-y_t)^\trns Q_k(x_t^i-y_t)$, $\Vert u_t^i\Vert^2_{R_k}=u_t^{i\trns}R_k u_t$, $\rho>0$ is a scalar discount factor, and all matrix and vector coefficients are of appropriate dimension.
\begin{assum}[Convexity]\label{convexityCondMinorLQGMFG}
$R_k >0$, and $Q_k - S_k R_k^{-1} S_k^\intercal\geq 0$, for $k \in \mfK  $.
\end{assum}

From \eqref{MinorDyn} and \eqref{minorCost}, each agent interacts with all agents through the average state $x^{(N)}_t$ and the average control action $u^{(N)}_t$. For applications of such systems we refer the reader to \cite{FirooziISDG2017, casgrain2018mean, CarmonaSysRisk2015, JaimungalSysRisk2017} and the references therein. 

\subsubsection*{Control $\sigma$-Fields and Admissible Control Sets}

We denote the natural filtration generated by agent-$i$'s Wiener process $(w^i_{t})_{t\in\mfT}$ and initial state $x^i_0$ by $\mc{F}^i\coloneqq (\mc{F}_{t}^i)_{t\in\mfT}$, $ i \in \mfN$. 
 
\begin{assum}[Admissible Control Set] \label{ass: MinorContrAction}
 For each minor agent-$i,\, i\in \mfN$, the set of control inputs $\mc{U}^{i}$ is defined to be the collection of control laws $u^{i}_t$ adapted to the filtration
 $\mc{F}^{i} \coloneqq (\mc{F}_{t}^{i})_{t\in\mfT}, \, \mfT = [0, \infty),$ such that  
 \begin{itemize}
\item[(i)] $\mb{E}[\int_0^\infty e^{-\rho t} u_t^{i\intercal}  u_t^i\, dt] < \infty$,
\item[(ii)] $x^i_t$ satisfies the dynamics \eqref{MinorDyn}, and the cost functional \eqref{minorCost} is finite.
 \end{itemize}
\end{assum}
\subsection{Infinite Populations}
In the mean field game methodology, the optimization problem \eqref{MinorDyn}-\eqref{minorCost} is solved in the infinite population limit where the number of agents $N$ goes to infinity. This results in a simpler stochastic optimal control problem for the generic agent-$i$ given by
\begin{align} \label{MinorDyn_inf}
&dx_t^i = (A_k x_t^i + \bar{F}_k \bar{x}_t+ \bar{H}_k\bar{u}_t + B_k u_t^i + b_k(t))dt + D_k dw_t^i,\\
&J^{\infty}_i(u^i) = \mb{E} \bigg[ \int_0^{\infty} e^{-\rho t}\bigg( \hf \Vert x_t^i-\bar{y}_t\Vert^2_{Q_k} + \hf \Vert u_t^i\Vert^2_{R_k}\nonumber\\&\hspace{1.5cm}+ (x_t^i-\bar{y}_t)^\intercal S_k u_t^i + \eta_k^\trns (x_t^i-\bar{y}_t) + n_k^\trns u_t^{i} \bigg) dt \bigg],\label{Minorcost_inf}
\end{align}
where $\bar{y}_t = \bar{\psi_k} \bar{x}_t$, $\bar{F}_k:= F_k \otimes \begin{bmatrix}\pi_1,...,\pi_K\end{bmatrix}$, $\bar{H}_k:= H_k \otimes \begin{bmatrix}\pi_1,...,\pi_K\end{bmatrix}$, $\bar{\psi}_k:= \psi_k \otimes \begin{bmatrix}\pi_1,...,\pi_K\end{bmatrix}$. Moreover, $(\bar{x}_t)^\trns = \begin{bmatrix} (\bar{x}_t^{1})^{\trns}, \dots, (\bar{x}_t^{K})^{\trns}\end{bmatrix}$, $\bar{u}_t^\trns = \begin{bmatrix} (\bar{u}_t^1)^{\trns}, \dots, (\bar{u}_t^K)^\trns\end{bmatrix}$ with $\bar{x}_t^k := \lim_{N_k\rightarrow \infty} \frac{1}{N_k}\sum_{i\in \mc{I}_k} x^{i,k}_t$, and $\bar{u}_t^{k} := \lim_{N_k\rightarrow \infty} \frac{1}{N_k}\sum_{i\in \mc{I}_k} u^{i,k}_t$, given that the limits exist. From \eqref{MinorDyn_inf}, all agents are connected through $\bar{x}_t$ and $\bar{u}_t$, which can be deterministically specified by a generic agent. The solutions to the limiting problem and their connection with the finite-population case are described in the next section.
\subsection{Mean Field Consistency and $\epsilon$-Nash Equilibrium}
We denote by $\mathbf{e}_k \in \dsR^{n \times nK}$ the matrix which has the identity matrix $I_n$ in the $k$th block and zero matrix $0_{n \times n}$ in the other $(K-1)$ blocks, i.e.
$\mathbf{e}_k = \left[ 0_{n\times n}, ... ,  0_{n\times n}, I_n,  0_{n\times n}, ...,  0_{n\times n}\right ]$. We also define 
\begin{align}
    \bar{A} = \begin{bmatrix}
    \bar{A}_1\\
    \vdots\\
    \bar{A}_K
    \end{bmatrix},
\end{align}
where $\bar{A}_k,\, k \in \mfK$ is given by \eqref{eq:barA_k}.
\begin{assum} \label{ass:MFEquationSolAss}
The parameters in \eqref{MinorDyn}-\eqref{minorCost} are such that 
\begin{itemize}
    \item[(i)] the resulting set of mean field consistency equations
\eqref{fixed_point_eqs_classic} have a unique bounded solution ($\Pi_k$, $s_k$, $\bar{A}_k$, $\bar{m}_k$), and
\item[(ii)] we have \begin{gather*}
\Pi_k > 0,\quad \bar{A} - \tfrac{\rho}{2} I_{nK} <0,
\\
A_k- B_kR_k^{-1}S_k^\trns-B_k R_k^{-1}B_k^\trns \Pi_k -\tfrac{\rho}{2} I_n <0.
\end{gather*}
\end{itemize}
 \end{assum}
As we have an infinite-horizon problem, the second part of this assumption is necessary to ensure the cost functional is finite. 

 
The following theorem encapsulates the classical results for a general case of LQG MFG systems which have been studied extensively in the literature. It characterizes a set of best-response control actions that forms a Nash-equilibrium and an $\epsilon$-Nash equilibrium for, respectively, the infinite-population and the finite-population LQG MFG system. The equilibrium is tied to the fixed-point solution to the mean-field consistency equations determining the values of the parameters used in the control actions.
\begin{thm}[Classical LQG MFG Solutions]\label{thm:ClassicLQGMFG}
Subject to \textit{Assumptions \ref{IntialStateAss}-\ref{ass:MFEquationSolAss}}, the system of equations \eqref{MinorDyn_inf}-\eqref{Minorcost_inf} together with the consistency equations
\begin{subequations}\label{fixed_point_eqs_classic}
\begin{align}
\rho &\Pi_k = \Pi_k A_k + A_k^\trns \Pi_k \nonumber\allowdisplaybreaks\\ 
& \hspace{1.5cm}-(\Pi_k B_k+S_k)R_k^{-1}(B_k^\trns \Pi_k +S_k^\trns) +Q_k,
 \allowdisplaybreaks\\
\rho &s_k(t) = \dot{s}_k(t)+(A_k^\trns- S_k R_k^{-1}B_k^\trns-\Pi_k B_k R_k^{-1}B_k^\trns)s_k(t) 
\nonumber\allowdisplaybreaks\\
&\hspace{0.2cm} +\Pi_k\big((\bar{F}_k\!+\!B_k R_k^{-1} S_k^\trns \bar{\psi})\bar{x}_t\!\!+\!\bar{H}_k \bar{u}_t \!-\!B_kR_k^{-1}n_k\!+\! b_k(t)\big)\nonumber \allowdisplaybreaks\\&\hspace{1.1cm}
+(S_kR_k^{-1}S_k^\trns- Q_k)\bar{\psi}_k\bar{x}_t - S_k R_k^{-1}n_k+\eta_k,
\label{offset-exp}\allowdisplaybreaks\\
&\bar{A}_k = \left( A_k- B_k R_k^{-1}(B_k^\trns \Pi_k + S_k^\trns)\right)\mathbf{e}_k \nonumber\\&\hspace{3cm} +B_k R_k^{-1}S_k^\trns \bar{\psi}_k + \bar{F}_k +  \bar{H}_k J,\label{eq:barA_k}\allowdisplaybreaks\\
&\bar{m}_k = - B_k R_k^{-1}(B_k^\trns s_k(t)+n_k) +\bar{H}_k L +b_k(t),\allowdisplaybreaks\\
&d\bar{x}_t^{k} = (\bar{A}_k \bar{x}_t^{k} + \bar{m}_k)dt,\allowdisplaybreaks\\
& \bar{u}_t = J \bar{x}_t + L,\label{ubar_classic}\allowdisplaybreaks\\
& J = \begin{bmatrix}
-R_1^{-1}(B_1^\trns \Pi_1 + S_1^\trns)\mathbf{e}_1 + R_1^{-1}S_1^\trns\bar{\psi}_1  \\
\vdots\\
-R_K^{-1}(B_K^\trns \Pi_K + S_K^\trns)\mathbf{e}_K + R_K^{-1}S_K^\trns\bar{\psi}_K 
\end{bmatrix}, \label{J}\allowdisplaybreaks\\
&L = \begin{bmatrix}
-R_1^{-1}(B_1^\trns s_1(t)+n_1) \\
\vdots\\
-R_K^{-1}(B_K^\trns s_K(t)+n_K)
\end{bmatrix},\label{L}
\end{align}
\end{subequations}
generate a set of actions $\mc{U}^{\infty}_{MF}:=\!\{u^{i,\ast}\!\!, i=1,\dots, \infty\}$ given by  \begin{equation}\label{opt-mean_classic}
u_t^{i,\ast} = -R_k^{-1}\left[(B_k^\trns \Pi_k + S_k^\trns) x_t^{i} + (B_k^\trns s_k(t)-S_k^\trns\bar{\psi}_k \bar{x}_t+n_k)\right],
\end{equation}
such that
\begin{enumerate}
\item[(i)]  The set of infinite population control laws $\mc{U}^{\infty}_{MF}$ yields the infinite population Nash equilibrium, i.e.,
\begin{equation*}
 J_i^{\infty}(u^{i,*}, u^{-i,*}) = \inf_{u^i \in \mc{U}^{i}}J_i^{\infty} (u^i, u^{-i,*}).
  \end{equation*}
\item[(iii)] The set of control laws $\mc{U}_{MF}^{N} := \{ {u}^{i,*};  i \in \mfN \}$, $1\leq N < \infty$, forms an $\epsilon$-Nash equilibrium for all $\epsilon$, i.e., for all $\epsilon>0$, there exists $N(\epsilon)$ such that for all $N \geq N(\epsilon)$ 
\begin{equation*}
J_i^{N}(u^{i,*},  u^{-i,*})-\epsilon \leq\inf_{u^i \in\mc{U}^i} J_i^{N}(u^i, u^{-i,*}) \leq  J_i^{N}( u^{i,*}, u^{-i,*}).
\end{equation*}
$\hfill \square$
\end{enumerate}
\end{thm}
\textit{Proof.}
 Following the approach in \cite{HuangTAC2007}, the results could be generalized to the systems given by \eqref{MinorDyn}-\eqref{minorCost}. Moreover, the system \eqref{MinorDyn}-\eqref{minorCost} could be considered as a special case of the major-minor LQG MFGs studied in \cite{FCJ-Convex2018} (for a comparison of different approaches to major-minor LQG MFGs see \cite{Firoozi2020,HuangCIS2020}.).   $\hfill \square$
 
\begin{remark} If \textit{Assumption \ref{ass:MFEquationSolAss}} holds for $\rho=0$, for agent-$i$ governed by \eqref{MinorDyn_inf}-\eqref{Minorcost_inf}, we have $\mb{E}\left[\int_0^t (x_t^{i})^\trns x_t^{i}\right]<\infty$ and $\mb{E}\left[\int_0^t (u_t^{i,\ast})^\trns u_t^{i,\ast}\right]<\infty$ for each $t \in \mfT$. 
\end{remark}
\section{Exploratory LQG MFGs with Entropy Regularization}\label{sec:reg_LQGMFG}

Following the approach in \cite{XYZRL2020}, we extend the notion of action to distribution over actions. Consider $n$ identical and independent rounds of experiments for agent $i$, where at round $j,\, j \in \{1,\dots,n\}$, the action $u^{i,j}$ is sampled with respect to the distribution $\Phi^i_t(u)$ belonging to the space $\mathcal{P}(\dsR^m)$ of probability measures on $\dsR^m$, and is executed for the system \eqref{MinorDyn}-\eqref{minorCost}. The corresponding states are denoted by $x^{i,j}_t, \,j \in \{1, \dots, n \}$. The increments of the state over the interval $[t, t+\Delta t]$ are given by
\begin{gather*}
\Delta x_t^{i,j} = x_{t+\Delta t}^{i,j} - x_t^{i,j}. 
\end{gather*}
As $n \rightarrow \infty$, by the law of large numbers we have
\begin{gather*}
\frac{1}{n} \sum_{j=1}^{n} \Delta x_t^{i,j} \rightarrow \mathbb{E}\left[\Delta x_t^{i,\Phi}\right],\quad a.s.\\
\frac{1}{n} \sum_{j=1}^{n} u^{i,j}_t \,\Delta t \rightarrow \mathbb{E}\left[\int_{\mathds{R}^m}u\; \Phi_t^i(u)\;du \right] \Delta t = \mathbb{E}\left[ \mu_t^i \right] \Delta t, \quad a.s.
\end{gather*}
where $x_t^{i,\Phi}$ denotes the resulting state when the agent-$i$'s actions are sampled from the distribution $\Phi^i_t$, and $\mu_t^i := \int_{\mathds{R}^m} u\; \Phi_t^i(u)\,du$. Moreover, for every term in the objective function $r(x^{i,j}_t, u_t^{i,j})$, we can write 
\begin{align*}
\frac{1}{n}\, \sum_{j=1}^{n}r(x^{i,j}_t, u_t^{i,j}) \Delta t \rightarrow \mathbb{E}\left[\int_{\mathds{R}^m} r(x^{i,\Phi}_t, u)\Phi_t^i(u)du \right] \Delta t, \quad a.s.
\end{align*}
We use the above exposition to derive the exploratory dynamics and cost functional for agent-$i$, $i\in \mfN$. 
\subsection{Finite Populations}\label{sec:finPopExp}
We introduce the exploratory dynamics for agent-$i$ in the finite population as in 
\begin{multline}\label{dyn_modified_fin}
dx_t^{i,\Phi^N} = \Big(A_k x_t^{i,\Phi^N} + F_k x_t^{(N),\Phi^N} + H_k \mu^{(N)}_t \\+ B_k \mu_t^{i,N} + b_k(t)\Big)dt + D_k dw_t^i,
\end{multline}
$i \in \mfN$, $k \in \mfK$, where $x_t^{i,\Phi^N} \in \dsR^n$ denotes the resulting state when actions are sampled with respect to the distribution $\Phi_t^{i,N} \in \mathcal{P}(\dsR^m)$ and executed. Moreover, $\mu_t^{i,N} = \int_{\mathds{R}^m} u \Phi_t^{i,N}(u)du$,
$x_t^{(N),\Phi^N} = \frac{1}{N} \sum_{i\in \mfN} x^{i,\Phi^N}_t$, and
$\mu_t^{(N)} = \frac{1}{N} \sum_{i \in \mfN} \mu^{i,N}_t$.

The exploratory cost functional for agent-$i$ is introduced as 
\begin{multline}\label{cost_exp_fin_pop}
J_i(\Phi^{i,N},\Phi^{-i,N}) = \mb{E} \bigg[\int_0^{\infty} e^{-\rho t}\bigg(\hf\left\Vert x_t^{i,\Phi^N}- y_t^{\Phi^N}\right\Vert^2_{Q_k}\allowdisplaybreaks\\ +\eta_k^\trns (x_t^{i,\Phi^N}-y_t^{\Phi^N})+ \phi_k \Big(\int_{\dsR^m} \Phi_t^{i,N}(u)du -1\Big)\allowdisplaybreaks\\+\int_{\dsR^m}\Big((x_t^{i,\Phi^N}-y_t^{\Phi^N})^\trns S_k + \hf u^\trns R_k + n_k^\trns\Big) u\Phi_t^{i,N}(u)du\bigg) dt \bigg],
\end{multline}
where 
\begin{equation}
y_t^{\Phi^N}=\psi_k x_t^{(N),\Phi^N}. 
\end{equation}
The third term in \eqref{cost_exp_fin_pop} is incorporated to enforce the integral constraint $\int_{\dsR^m}\Phi^{i,N}_t(u) du=1, \forall t \in \mfT$, as $\phi_k \rightarrow \infty$, where $\phi_k$ is a Lagrange multiplier. All other processes and the coefficients in \eqref{dyn_modified_fin}-\eqref{cost_exp_fin_pop} are as defined in \eqref{MinorDyn}-\eqref{minorCost}.
 
To reward exploration we incorporate Shannon's differential entropy given by
\begin{align*}
\mathcal{H}(\Phi^{i,N}) = - \int_{\dsR^m} {\Phi^{i,N}(u)}\ln {\Phi^{i,N}(u)}du,
\end{align*}
in the cost functional, which measures the level of exploration. Hence the exploratory entropy-regularized cost functional for agent-$i$ is given by 
\begin{multline}\label{cost_exp_fin_pop_reg}
J_i(\Phi^{i,N},\Phi^{-i,N}) = \mb{E} \bigg[\int_0^{\infty} e^{-\rho t}\bigg(\hf\left\Vert x_t^{i,\Phi^N}- y_t^{\Phi^N}\right\Vert^2_{Q_k}\allowdisplaybreaks\\ +\eta_k^\trns (x_t^{i,\Phi^N}-y_t^{\Phi^N})+ \phi_k \Big(\int_{\dsR^m} \Phi_t^{i,N}(u)du -1\Big)\allowdisplaybreaks\\+\int_{\dsR^m}\Big(\big((x_t^{i,\Phi^N}-y_t^{\Phi^N})^\trns S_k + \hf u^\trns R_k + n_k^\trns\big)u \allowdisplaybreaks\\+ \lambda_k \ln \Phi_t^{i,N}(u)\Big) \Phi_t^{i,N}(u)du\bigg) dt \bigg],
\end{multline}
where $\lambda_k$ represents the weight associated with exploration. 

\subsubsection*{Admissible Control Set}
For the exploratory LQG MFGs, we define the set of admissible control distributions $\mc{U}^{i,\tex}$ to include $\Phi^i_t \in \mathcal{P}(\mathds{R}^m)$ such that 
\begin{itemize}
    \item for each $A \in \mathcal{B}(\dsR^m)$, $\{\int_{A} \Phi^{i}_t(u)du, t \in \mfT \}$ is $\mc{F}_t^i$-measurable;
    \item $\mb{E}\left[\int_0^\infty e^{-\rho t}\mu_t^{i\intercal}\mu_t^{i} dt\right]< \infty$, where $\mu_t^{i} = \int_{\mathds{R}^m} u \Phi_t^{i}(u)du$;
    \item for each $t \in \mfT$, $x^{i,\Phi}_t$ satisfies \eqref{dyn_modified_fin}, and the cost functional \eqref{cost_exp_fin_pop_reg} is finite.
\end{itemize}
\begin{obs} Only the mean value $\mu^{i,N}_t$ of the distribution $\Phi_t^{i,N}$ appears in the dynamics \eqref{dyn_modified_fin}. Intuitively, the dynamics \eqref{dyn_modified_fin} describe the average behavior of a single agent governed by the dynamics \eqref{MinorDyn} who samples actions according to $\Phi_t^{i,N}$ and applies them to \eqref{MinorDyn}. This is to be distinguished from the mass behavior of all agents which is described by the mean field distribution in the infinite-population case where the number of agents $N$ goes to infinity. 
\end{obs}
We aim to find the optimal control distribution under which an agent samples and applies random control actions. The execution of random actions is usually used for learning purposes. In this case it could be used to learn unknown parameters of the dynamics. Finding the optimal control distribution, which provides optimally randomized control actions, helps reduce the number of interactions with the system and could expedite the learning process. However, it is very difficult to find the optimal control distribution in the finite-population case due to the coupling of the agents with each other through the terms $x_t^{(N),\Phi^N}$ and $\mu_t^{(N)}$. Hence we first solve the problem in the infinite-population limit and then relate the resulting solutions to the finite-population case.      
\subsection{Infinite Populations}
In the infinite population limit, where the number of agents goes to infinity ($N \rightarrow \infty$), the exploratory dynamics for agent-$i$ are given by
\begin{align}\label{dyn_modified_inf}
dx_t^{i,\Phi} = \left(A_k x_t^{i,\Phi} +\bar{F}_k \bar{x}_t^{\Phi} + \bar{H}_k \bar{\mu}_t + B_k \mu_t^i + b_k(t)\right)dt + D_k dw_t^i,
\end{align}
where $x_t^{i,\Phi} \in \dsR^n$ denotes the resulting state when actions are sampled with respect to the distribution $\Phi_t^{i} \in \mathcal{P}(\dsR^m)$ and executed, $\mu^i_t = \int_{\mathds{R}^m} u \Phi_t^{i}(u)du$, $(\bar{x}_t^{\Phi})^\trns = \begin{bmatrix} (\bar{x}_t^{1,\Phi})^{\trns}, \dots, (\bar{x}_t^{K,\Phi})^{\trns}\end{bmatrix}$, $\bar{\mu}_t^\trns = \begin{bmatrix} (\bar{\mu}_t^1)^{\trns}, \dots, (\bar{\mu}_t^K)^\trns\end{bmatrix}$, with $\bar{x}_t^{k,\Phi}\!\! :=\! \lim_{N_k \rightarrow \infty} \frac{1}{N_k} \sum_{i\in \mc{I}_k} x^{i,\Phi}_t$ and $\bar{\mu}_t^k\! :=\! \lim_{N_k \rightarrow \infty} \frac{1}{N_k} \sum_{i \in \mc{I}_k} \mu^i_t$ given that the limits exist. 

The entropy-regularized exploratory cost functional for agent-$i$ in the infinite-population limit is given by
\begin{multline}\label{cost_func_exploratory}
J_i^\infty(\Phi^i) = \mb{E} \bigg[ \int_0^{\infty} e^{-\rho t}\bigg(\hf \Vert x_t^{i,\Phi}- \bar{y}^{\phi}_t\Vert^2_{Q_k} + \eta_k^\trns (x_t^{i,\Phi}-\bar{y}_t^{\Phi})\allowdisplaybreaks\\ + \phi_k \Big(\int_{\dsR^m} \Phi_t^i(u)du -1\Big) + \int_{\dsR^m}\Big( \big(({x_t^{i,\Phi}}-\bar{y}_t^{\Phi})^\trns S_k \allowdisplaybreaks\\+ \hf u^\trns R_k  + n_k^\trns \big)u + \lambda_k \ln \Phi_t^i(u) \Big)\Phi_t^i(u)du \bigg) dt \bigg],
\end{multline}
where $\bar{y}_t^{\Phi} = \bar{\psi}_k \bar{x}_t^{\Phi}$. All other processes and the coefficients in \eqref{dyn_modified_inf}-\eqref{cost_func_exploratory} are as defined in \eqref{MinorDyn_inf}-\eqref{Minorcost_inf}.

To obtain the optimal control distribution for agent-$i$, we first compute the Gâteaux derivative of the cost functional \eqref{cost_func_exploratory}.
\begin{thm}[G\^ateaux Derivative] The G\^ateaux derivative of agent-$i$'s cost functional \eqref{cost_func_exploratory} in an arbitrary direction $\omega^i$, such that the perturbed distribution $\Phi^{i,\eps}_t = e^{\eps\omega^i_t(u)}\Phi^i_t(u) \in \mc{U}^{i,\tex}$, is given by 
\begin{align}\label{gat-deriv-reduced}
\langle \mc{D}J_i^\infty(\Phi^i),\omega^i \rangle &=\lim_{\eps\rightarrow 0} \frac{J_i^\infty(\Phi^{i,\eps})-J_i^\infty(\Phi^{i})}{\eps}\nonumber\\
=\mb{E}\bigg[\int_0^{\infty}& e^{-\rho t} \bigg( \int_{\dsR^m}  f(t, u, x^{i,\Phi}_t, \Phi^i_t)\Phi^i_t(u)\omega_t^i(u)du \bigg) dt \bigg], 
\end{align}
where 
\begin{multline}
f(t,u,x^{i,\Phi}_t,\Phi^i_t)=\bigg(\mb{E}\Big[\int_t^\infty e^{-\rho (s-t)}  \big( Q_k(x_s^{i,\Phi} - \bar{y}_s^\Phi) \allowdisplaybreaks\\+ \eta_k + S_k \mu_s^i\big)^\intercal e^{A_k(s-t)}ds \Big\vert \mc{F}^i_t \Big]  B_k   +(x_t^{i,\Phi}-\bar{y}_t^{\Phi})^\trns S_k \allowdisplaybreaks\\+ \hf u^\trns R_k +n_k^\trns\bigg) u + \lambda_k \ln\Phi_t^i(u)+\lambda_k + \phi_k.
\end{multline}
$\hfill \square$
\end{thm}
\textit{Proof.} We define 
\begin{equation}
\Phi_t^i(u) = e^{g_t^i(u)}, 
\end{equation}
and rewrite the cost functional as
\begin{multline}
J_i^\infty(\Phi^i) = \mb{E} \bigg[ \int_0^{\infty} e^{-\rho t}\bigg(\hf \Vert x_t^{i,\Phi}-\bar{y}_t^{\Phi}\Vert^\trns_{Q_k} + \eta_k^\trns (x_t^{i,\Phi}-\bar{y}_t^{\Phi}) \allowdisplaybreaks\\
+ \int_{\dsR^m}\Big(\big((x_t^{i,\Phi}-\bar{y}_t^{\Phi})^\trns S_k + \hf u^\trns R_k 
+ n_k^\trns\big) u \allowdisplaybreaks\\+ \lambda g_t^i(u) \Big)e^{g_t^i(u)}du 
+ \phi_k \Big(\int_{\dsR^m} e^{g_t^i(u)}du -1\Big)\bigg) dt \bigg]. 
\end{multline}
SDE \eqref{dyn_modified_inf} may be re-written as
\begin{multline}\label{tsol_exp}
x_t^{i,\Phi} = e^{A_kt} x_0^i + \int_0^t e^{A_k(t-s)} (\bar{F}_k \bar{x}_s^{\Phi} + \bar{H}_k \bar{\mu}_s +b_k(s) )ds\\+ \int_0^t e^{A_k(t-s)}B_k \mu_s^i ds + \int_0^t e^{A_k(t-s)}D_kdw_s^i.
\end{multline}
We perturb $g_t^i(u)$ in an arbitrary direction $\omega^i_t$ such that $e^{g^i_t+\epsilon\omega^i_t} \in \mc{U}^{i,\tex},\, \forall t \in \mfT$, as in 
\begin{equation}
g_t^{i,\epsilon}(u) = g_t^i(u) + \epsilon \omega_t^i(u),
\end{equation}
and use the Taylor expansion of $e^{\eps \omega_t^i(u)}$ to write 
\begin{align}
\mu_t^{i,\eps} &= \int_{\dsR^m} u e^{g_t^i(u) + \epsilon \omega_t^i(u)} du,\allowdisplaybreaks\\
&= \int_{\dsR^m} u e^{g_t^i(u)}(1 + \eps \omega_t^i(u)+O(\eps^2)) du,\allowdisplaybreaks\\
& = \mu_t^i + \eps \int_{\dsR^m} u e^{g_t^i(u)} \omega_t^i(u) du + O(\eps^2).
\end{align}
From \eqref{tsol_exp}, the corresponding perturbed state $x_t^{i,\Phi, \eps}$ satisfies
\begin{gather}
x_t^{i,\Phi, \eps} = x_t^{i,\Phi} + \eps \tilde{\omega}_{t}^i + O(\eps^2),\\
\tilde{\omega}_{t}^i \coloneqq  \int_0^t e^{A_k(t-s)} B_k \int_{\dsR^m} u e^{g_s^i(u)} \omega_s^i(u) du ds.\label{omega_tilde}
\end{gather}
%
%
The perturbed cost functional can then be written as  
\begin{multline}
J_i^{\infty}(\Phi^{i,\eps}) = J_i(\Phi^i) + \eps \mb{E} \bigg[ \int_0^{\infty} e^{-\rho t}\bigg(\Big(Q_k (x_t^{i,\Phi}-\bar{y}_t^\Phi) + \eta_k\allowdisplaybreaks\\ +S_k \mu^{i}_t \Big)^\intercal \tilde{\omega}_t^i
+ \int_{\dsR^m}\Big(\big((x_t^{i,\Phi}-\bar{y}_t^\Phi)^\trns S_k + \hf u^\trns R_k +n_k^\trns\big) u \allowdisplaybreaks\\+ \lambda_k g_t^i(u)+\lambda_k +\phi_k\Big) e^{g_t^i(u)} \omega_t^i(u)du \bigg) dt \bigg] + O(\eps^2). 
\end{multline}
By substituting \eqref{omega_tilde} and changing the order of integration, we have
\begin{multline}\label{term1_red}
\mb{E}\int_0^{\infty} e^{-\rho t}\Big(Q_k(x_t^{i,\Phi}-\bar{y}_t^\Phi) + \eta_k + S_k\mu^{i}_t  \Big)^\intercal\tilde{\omega}_t^i dt \allowdisplaybreaks\\=\mb{E}\int_0^{\infty} e^{-\rho t} \bigg( \int_{\dsR^m} \mb{E}\Big[\int_t^\infty  e^{-\rho (s-t)} 
\Big(Q_k(x_s^{i,\Phi}-\bar{y}_s^\Phi) + \eta_k \allowdisplaybreaks \\+ S_k\mu_s^i\Big)^\trns e^{A_k(s-t)}ds \Big| \mc{F}^i_t\Big] B_k u e^{g_t^i(u)} \omega_t^i(u)du\bigg)dt.
\end{multline}
Hence, \eqref{gat-deriv-reduced} is proven. $\hfill \square$

Now we use \eqref{gat-deriv-reduced} to compute the optimal control distribution for agent-$i$ as stated below. 
\begin{thm}[Optimal Control Distribution]\label{thm:opt_dist} Subject to \textit{Assumptions \ref{IntialStateAss}-\ref{ass:MFEquationSolAss}}, the optimal control distribution for agent-$i$, with dynamics \eqref{dyn_modified_inf} and cost functional \eqref{cost_func_exploratory}, is given by 
\begin{align}\label{opt-dist}
\Phi_t^{i,\ast}(u) = \frac{\exp\left[-\frac{1}{2\lambda_k}(u-\mu_t^{i,\ast})^\trns R_k (u-\mu_t^{i,\ast})\right]}{\int_{\mathds{R}^m}\exp\left[-\frac{1}{2\lambda_k}(u-\mu_t^{i,\ast})^\trns R_k (u-\mu_t^{i,\ast})\right]du},
\end{align}
where
\begin{equation}\label{opt-mean}
\mu_t^{i,\ast} = -R_k^{-1}\left[(B_k^\trns \Pi_k + S_k^\trns) x_t^{i,\Phi} + (B_k^\trns s_k(t)-S_k^\trns\bar{\psi}_k \bar{x}_t^{\Phi}+n_k)\right].
\end{equation}
Furthermore, the matrices $\Pi_k$ and $s_k(t)$ are the fixed point solutions to the following set of equations.
\begin{subequations}\label{fixed_point_eqs}
\begin{align}
\rho &\Pi_k = \Pi_k A_k + A_k^\trns \Pi_k -(\Pi_k B_k+S_k)\nonumber \\
&\hspace{3cm}\times R_k^{-1}(B_k^\trns \Pi_k+S_k^\trns) +Q_k,
 \allowdisplaybreaks\\
\rho &s_k(t) = \dot{s}_k(t)+(A_k^\trns- S_k R_k^{-1}B_k^\trns-\Pi_k B_k R_k^{-1}B_k^\trns)s_k(t) 
\nonumber\allowdisplaybreaks\\
&\hspace{0.2cm} +\Pi_k\big((\bar{F}_k\!+\!B_k R_k^{-1} S_k^\trns \bar{\psi})\bar{x}_t^{\Phi}\!\!+\!\bar{H}_k \bar{\mu}_t \!-\!B_kR_k^{-1}n_k\!+\! b_k(t)\big)\nonumber \allowdisplaybreaks\\&\hspace{1.1cm}
+(S_kR_k^{-1}S_k^\trns- Q_k)\bar{\psi}_k\bar{x}_t^\Phi - S_k R_k^{-1}n_k+\eta_k,
\label{offset-exp}\allowdisplaybreaks\\
&\bar{A}_k = \left( A_k- B_k R_k^{-1}(B_k^\trns \Pi_k + S_k^\trns)\right)\mathbf{e}_k \nonumber\\&\hspace{3cm} +B_k R_k^{-1}S_k^\trns \bar{\psi}_k + \bar{F}_k +  \bar{H}_k J,\allowdisplaybreaks\\
&\bar{m}_k = - B_k R_k^{-1}(B_k^\trns s_k(t)+n_k) + b_k(t)+  \bar{H}_k L,\allowdisplaybreaks\\
&d\bar{x}_t^{k,\Phi} = (\bar{A}_k \bar{x}_t^{k,\Phi} + \bar{m}_k)dt, \label{stateMF_mean_exp}\allowdisplaybreaks\\
&\bar{\mu}_t = J \bar{x}_t^{\Phi} + L \label{cntlMF_mean_exp},
\end{align}
where $J$ and $L$ are given by \eqref{J} and \eqref{L}. 
\end{subequations}
$\hfill \square$
\end{thm}
\textit{Proof.} We note that the cost functional \eqref{cost_func_exploratory} is strictly convex with respect to the distribution $\Phi^i_t(u)$. Hence the necessary and sufficient condition for $\Phi^{i,\ast}_t(u)$ to be the optimal distribution for agent-$i$ is given by
\begin{equation*}
\langle \mc{D}J_i^\infty(\Phi^{i,\ast}),\omega^i \rangle =0, 
\end{equation*}
for every $\omega^i$ such that $\Phi^{i,\ast,\eps}_t = e^{\eps\omega^i_t(u)}\Phi^{i,\ast}_t(u) \in \mc{U}^{i,\tex}$. For the G\^ateaux derivative \eqref{gat-deriv-reduced} to be zero for every admissible process $\omega^i_t$, we must have 
\begin{equation}
f(t,u,x^{i,\Phi}_t,\Phi^{i,\ast}_t) = 0.
\end{equation}
This results in
\begin{align}\label{exp_dist}
\Phi_t^{i,\ast}(u) &= \dfrac{\exp\left[- \frac{1}{2\lambda_k}\left(u - \mu^{i,\ast}_t \right)^\trns R_k \left(u - \mu^{i,\ast}_t\right)\right]}{\int_{\dsR^m} \exp\Big[- \frac{1}{2\lambda_k}\left(u- \mu^{i,\ast}_t \right)^\trns R_k \left( u - \mu^{i,\ast}_t \right)\Big]du},
\end{align}
where 
\begin{gather}
\mu^{i,\ast}_t := - R_k^{-1} \left(B_k^\trns p^i_t + S_k^\trns (x^{i,\Phi}_t -\bar{y}^{\Phi}_t)+ n_k\right),\allowdisplaybreaks\\
p_t^i \coloneqq \mb{E}\Big[\int_t^\infty  e^{-\rho (s-t)} e^{A_k^\trns(s-t)} \tilde{p}^i_s\, ds \Big\vert \mc{F}^i_t\Big],\allowdisplaybreaks\\
 \tilde{p}^i_s := Q_k (x_s^{i,\Phi}-\bar{y}^{\Phi}_s) + \eta_k +  S_k \mu_s^{i,\ast}.
\end{gather}
We write $p_t^i$ in terms of the martingale $M^i_t$ as in 
\begin{align}
p^i_t &=  e^{\rho t}e^{- A_k^\trns t}\Big(M_t^i - \int_0^t  e^{-\rho s}e^{A_k^\trns s} \tilde{p}^i_s ds\Big),\label{adjoint}
\allowdisplaybreaks \\
&M_t^i \coloneqq \mb{E}\left[\int_0^\infty  e^{-\rho s}e^{A_k^\trns s}\tilde{p}^i_s\, ds\Big \vert \mc{F}^i_t \right].
\end{align}
Next we adopt an ansatz for $p_t^i$ given by
\begin{equation}\label{ansatz}
p_t^i = \Pi_k(t) x_t^{i,\Phi} + s_k(t). 
\end{equation}
Using the above ansatz, the mean value of $\Phi_t^{i,\ast}$ is given by
\begin{equation}\label{cntrlAve}
\mu_t^{i,\ast} = -R_k^{-1}\left[(B_k^\trns \Pi_k(t) + S_k^\trns) x_t^{i,\Phi} + (B_k^\trns s_k(t)-S_k^\trns \bar{y}^{\Phi}_t+n_k)\right].
\end{equation}
Now we verify that \eqref{exp_dist} belongs to the admissible control set specified in Section \ref{sec:finPopExp}. By inspection, \eqref{exp_dist} and its integral over every $A \in \mathcal{B}(\dsR^m)$ are $\mc{F}_t^i$-measurable for every $t \in \mfT$. Given \textit{Assumption \ref{ass:MFEquationSolAss}}, by using the triangle inequality we can easily show that under the control distribution \eqref{exp_dist}, we have $\mb{E}\left[\int_0^\infty e^{-\rho t} (x_t^{i,\Phi})^\trns x_t^{i,\Phi}\right]<\infty$ and $\mb{E}\left[\int_0^\infty e^{-\rho t}(\mu_t^{i,\ast})^\trns \mu_t^{i,\ast}\right]<\infty$. Moreover, we can show that 
\begin{align}\label{entropyEval}
&\mb{E} \bigg[\int_0^{\infty} e^{-\rho t}\Big( \lambda_k \int_{\dsR^m} \Phi_t^{i,\ast}(u)\ln \Phi_t^{i,\ast}(u) du\Big) dt \bigg] \nonumber \allowdisplaybreaks\\&\hspace{1em}= \frac{\lambda_k}{2\rho}\ln(2\pi \lambda_k R^{-1}_k),
\end{align}
and
\begin{align}
\mb{E} &\bigg[\int_0^{\infty} e^{-\rho t}\int_{\dsR^m}\hf u^\trns R_k u \Phi_t^{i,\ast}(u)du dt \bigg]\nonumber\allowdisplaybreaks\\
&= \mb{E} \bigg[\int_0^{\infty} e^{-\rho t}\!\!\bigg(\!\int_{\dsR^m}\hf (u-\mu^{i,\ast}_t)^\trns R_k (u-\mu^{i,\ast}_t) \Phi_t^{i,\ast}(u)du  \nonumber \allowdisplaybreaks\\&\hspace{0.7cm}+ \hf{(\mu_t^{i,\ast})}^\trns R_k {\mu_t^{i,\ast}}\bigg) dt \bigg]\nonumber\allowdisplaybreaks\\
&= \tfrac{\lambda_k}{2\rho}+ \hf\mb{E} \bigg[\int_0^{\infty} e^{-\rho t}{(\mu_t^{i,\ast})}^\trns R_k {\mu_t^{i,\ast}}dt\bigg], \label{distVar}
\end{align}
where we have added and subtracted the term $\hf{(\mu_t^{i,\ast})}^\trns R_k {\mu_t^{i,\ast}}$ in \eqref{distVar}. Finally, using the triangle inequality and Cauchy-Schwarz inequality together with \eqref{cntrlAve}-\eqref{distVar}, we can show 
\begin{align*}
J_i(\Phi^{i}) &\leq \alpha(t)\mb{E} \bigg[\int_0^{\infty} e^{-\rho t} (x_t^{i,\Phi})^\trns x_t^{i,\Phi}  dt \bigg] + \beta(t)
< \infty, 
\end{align*}
where $\alpha(t)$ and $\beta(t)$ are bounded and deterministic functions of time.

We can show that $x_t^{i,\Phi}$ is normally distributed. Hence $\mu_t^{i,\ast}$ is also normally distributed. In the infinite-population limit, the distribution of the process $\mu^{i,\ast}_t$ for the generic agent-$i$, $i \in \mc{I}_k$, coincides with the control mean field distribution of subpopulation $k$, which is indeed the distribution of processes $\{\mu^i_t\}_{i \in \mc{I}_k}$ across all agents in subpopulation $k$. The mean value 
$\bar{\mu}_t^{k} = \lim_{N_k\rightarrow \infty}\tfrac{1}{N_k}\sum_{i \in \mc{I}_k}\mu^i_t$ of the control mean field distribution for subpopulation $k$ is given by 
\begin{equation}\label{cntrlAveMF}
\bar{\mu}_t^{k} = -R_k^{-1}\left[(B_k^\trns \Pi_k(t) + S_k^\trns) \bar{x}_t^{k,\Phi} + (B_k^\trns s_k(t)-S_k^\trns \bar{y}^{\Phi}_t+n_k)\right].
\end{equation}
Therefore we  have
\begin{equation}\label{cntrl_mf_exp}
\bar{\mu}_t = J \bar{x}_t^{\Phi} + L,
\end{equation}
where $\bar{\mu}_t^\intercal = \big[(\bar{\mu}_t^{1})^\intercal, \cdots, (\bar{\mu}_t^{K})^\intercal \big]$, and
\begin{gather}
J = \begin{bmatrix}
-R_1^{-1}(B_1^\trns \Pi_1 + S_1^\trns)\mathbf{e}_1 + R_1^{-1}S_1^\trns\bar{\psi}_1  \\
\vdots\\
-R_K^{-1}(B_K^\trns \Pi_K + S_K^\trns)\mathbf{e}_K + R_K^{-1}S_K^\trns\bar{\psi}_K 
\end{bmatrix}, \nonumber\allowdisplaybreaks\\ L= \begin{bmatrix}
-R_1^{-1}(B_1^\trns s_1(t)+n_1)\\
\vdots\\
-R_K^{-1}(B_K^\trns s_K(t)+n_K)
\end{bmatrix}.
\end{gather}
Next we aim to find the equations that $\Pi_k(t)$ and $s_k(t)$ satisfy. Using It\^{o}'s lemma and \eqref{dyn_modified_inf}, the ansatz \eqref{ansatz} satisfies the SDE
\begin{multline}
dp_t^i  = \big[d\Pi_k(t) + \Pi_k(t)A_k - \Pi_k(t) B_k R_k^{-1}B_k^\trns \Pi_k(t) \allowdisplaybreaks\\
- \Pi_k(t)B_kR_k^{-1}S_k^\trns \big] x^{i,\Phi}_t dt + \Big(\Pi_k(t)\bar{F}_k \bar{x}_t^\Phi+ \Pi_k(t) \bar{H}_k \bar{\mu_t} \allowdisplaybreaks\\
- \Pi_k(t) B_k R_k^{-1} (B_k^\trns s_k(t)-S_k^\trns \bar{y}_t^\Phi +n_k)  + \Pi_k b_k(t)\Big) dt + ds_k(t) \Big]\allowdisplaybreaks\\
 + \Pi_k(t) D_k dw_t^i. \label{diff1}
\end{multline}
From \eqref{adjoint}, \eqref{ansatz}, and the martingale representation theorem, the SDE that $p^i_t$ satisfies is given by
\begin{multline}
dp^i_t = \Big[\left(\rho I_n - A_k^\trns + S_k R_k^{-1}B_k^\trns \right)\Pi_k(t)-Q_k+S_kR_k^{-1}S_k^\trns \Big]\\\times x_t^{i,\Phi} dt+ \Big[\left(\rho I_n - A_k^\trns + S_k R_k^{-1}B_k^\trns\right)s_k(t)\allowdisplaybreaks\\+Q_k \bar{y}_t^\Phi - \eta_k-S_kR_k^{-1}S_k^\trns\bar{y}_t^\Phi+S_kR_k^{-1}n_k\Big]dt + q_k(t)dw_t^i, \label{diff2}
\end{multline}
where $I_n$ is the identity matrix of size $n$. Finally we equate \eqref{diff1} with \eqref{diff2} to get
\begin{align}
\rho \Pi_k(t) &= \dot{\Pi}_k(t) + \Pi_k(t)A_k + A_k^\trns \Pi_k(t) \nonumber\\&\hspace{0cm} -(\Pi_k(t)B_k+S_k)R_k^{-1}(B_k^\trns \Pi_k(t)+S_k^\trns) +Q_k,\label{RicTemp}
 \\
\rho s_k(t) &= \dot{s}_k(t)+\left( A_k^\trns- S_k R_k^{-1}B_k^\trns-\Pi_k B_k R_k^{-1}B_k^\trns \right)s_k(t)
 \nonumber
 \\
 &\hspace{-0.3cm}+ \Pi_k(t)\Big(\bar{F}_k \bar{x}_t^{\Phi} +\bar{H}_k \bar{\mu_t} +B_k R_k^{-1} S_k^\trns \bar{y}_t^\Phi -B_kR_k^{-1}n_k\nonumber\\&\hspace{-0.3cm}+ b_k(t)\Big) +(S_kR_k^{-1}S_k^\trns- Q_k) \bar{y}_t^\Phi - S_k R_k^{-1}n_k+\eta_k,
\end{align}
where $\bar{y}_t^{\Phi} = \bar{\psi}_k \bar{x}_t^{\Phi}$, and  $\dot{\Pi}_k(t)=0$ as agent-$i$ has an infinite-horizon discounted cost functional and all parameters in \eqref{RicTemp} are time-invariant. It remains to specify the mean value $\bar{x}_t^{k,\Phi}$ of the state mean field distribution for subpopulation $k$. 
By definition we have
\begin{align}
 \bar{x}_t^{k,\Phi} = \lim_{N_k \rightarrow \infty} \frac{1}{N_k}\sum_{i\in \mc{I}_k} x^{i, \Phi}_t.   
\end{align}
We first substitute \eqref{cntrlAve} in the minor agent's dynamics \eqref{dyn_modified_inf} to obtain the closed-loop dynamics for agent-$i$, $i \in \mc{I}_k$.
We then take the average over subpopulation $k$, and then the limit as $N\rightarrow \infty$, to get
\begin{equation}
d\bar{x}_t^{k,\Phi} = (\bar{A}_k \bar{x}_t^{k,\Phi} + \bar{m}_t)dt,
\end{equation}
where
\begin{align}
&\bar{A}_k = \left( A_k- B_k R_k^{-1}(B_k^\trns \Pi_k(t)+S_k^\trns)\right)\mathbf{e}_k \nonumber\\&\hspace{3.5cm}+B_k R_k^{-1}S_k^\trns \bar{\psi}_k + \bar{F}_k + \bar{H}_kJ,\\
&\bar{m}_k = - B_k R_k^{-1}(B_k^\trns s_k(t)+n_k) + \bar{H}_kL + b_k(t).
\end{align}
$\hfill \square$

\begin{remark} The results of Theorem \ref{thm:opt_dist} also hold when the dynamics of agent-$i$ given by \eqref{MinorDyn} are driven by a general martingale process defined on the same probability space where $w^i_t$ is defined. 
\end{remark}

\begin{remark} The results of Theorem \ref{thm:opt_dist} can be extended beyond the  Shannon entropy case. For example, consider replacing the Shannon entropy  by the relative entropy $\int_{\dsR^m}\Phi^i_t(u) \ln(\frac{\Phi^i_t(u)}{\Pi_t(u)})du$ -- which reduces to negative of Shannon entropy when the prior $\Pi_t(u)$ is uniform. With this change, the optimal exploratory control distribution becomes $\Pi_t(u)\Phi^{i,\ast}_t(u)$, where $\Phi^{i,\ast}_t(u)$ is given by \eqref{opt-dist}. In the interest of space, we skip the details of the computation. For a linear combination of Shannon entropy and cross entropy we refer the reader to \cite{Zariphopoulou2020}.   

\end{remark}

\begin{obs} We note that $\Phi^{i,\ast}_t$ is a Gaussian distribution with a stochastic mean value $\mu^{i,\ast}_t$ driven by the state sample path $x^i_t$ of the agent-$i$. From \eqref{opt-mean} and \eqref{opt-mean_classic}, $\mu^{i,\ast}_t$ coincides with the optimal control action $u^{i,\ast}_t$ of agent-$i$ in the classical setup \eqref{MinorDyn_inf}-\eqref{Minorcost_inf}. Sampling the actions with respect to $\Phi^{i,\ast}_t$ provides a randomization around the optimal control action of the agent-$i$ in the classical infinite-population limit setup. 
\end{obs}

\begin{obs} In the infinite-population limit, the law of state and the law of the control process that appears in the dynamics of a generic agent of subpopulation $k$ determine, respectively, the state mean field and the control mean field distributions of subpopulation $k$. The mixture of these distributions weighted by the proportion of agents in each subpopulation gives the mean field distributions across all agents. In the exploratory setup \eqref{dyn_modified_inf}-\eqref{cost_func_exploratory}, the law $\mc{L}(x^{.,k}_t)$ of the state of a represntative agent in subpopulation $k$ is Gaussian. Hence, the state mean field distribution $\sum_{k\in \mfK} \pi_k \mc{L}(x^{.,k}_t) $ is Gaussian with the mean value $\sum_{k \in \mfK}\pi_k\bar{x}_t^{k,\Phi}$, where $\bar{x}^{k,\Phi}_t$ is given by \eqref{stateMF_mean_exp}. Moreover, when a generic agent samples actions with respect to $\Phi^{i,\ast}_t$, according to the exploratory dynamics \eqref{dyn_modified_inf}, on average the mean value $\mu^{i,\ast}_t$ of $\Phi^{i,\ast}_t$ affects its dynamics. Hence the control mean field distribution in the exploratory case is specified by  $\sum_{k \in \mfK}\pi_k\mc{L}(\mu^{.,k,\ast}_t)$, which, from \eqref{opt-mean}, is Gaussian with the mean value $\sum_{k \in \mfK}\pi_k\bar{\mu}^{k,\ast}_t$, where $\bar{\mu}^{k,\ast}_t$ satisfies \eqref{cntlMF_mean_exp}. It is interesting to observe that both exploratory and classical LQG MFG systems share the same state and mean field distributions (this is formalized in Section \ref{sec:equivalence}). 


\end{obs}
\begin{obs} The exploratory control distribution $\Phi^{i,\ast}_t(u)$ of the generic agent-$i,\, i \in \mc{I}_k,$ is distinct from the control mean field distribution $\mc{L}(\mu^{.,k,\ast}_t)$ of subpopulation $k$. The former is a Gaussian distribution with the stochastic mean value $\mu^{i,\ast}_t$, while the latter is a Gaussian distribution with the deterministic mean value $\bar{\mu}^k_t$.

The mean $\mu^{i,\ast}_t$ of the exploratory control distribution $\Phi^{i,\ast}_t(u)$ is impacted by the mass effect of agents (i.e. the mean field) explicitly and through the offset term $s_k$ given by \eqref{offset-exp}. This offset term is a function the mean value $\sum_{k \in \mfK}\pi_k\bar{x}_t^{k,\Phi}$ of the state mean field. 


\end{obs}

In the next section we relate the solutions to the infinite-population case to that of the finite-population case. 
\subsection{$\epsilon$-Nash Property}
To show the $\epsilon$-Nash property for the set of obtained optimal control distributions $\{\Phi^{i,\ast}_t(u), i\in \mfN\}$, we first establish the relation between the infinite-population and finite-population state and cost functional for the generic agent-$i$.
\begin{thm} \label{thm:fin_inf_diff}Suppose that there exists a sequence $\{\delta_n \}_{n=1}^N$ such that $\delta_N \rightarrow 0$ as $N \rightarrow \infty$, and $\left| \tfrac{N_k}{N} -  \pi_k\right| = o(\delta_N)$, for all $k \in \mfK$. Given that all agents $j \in \mfN, j \neq i$, are using the optimal distributions with the corresponding mean values $\{\Phi^{j,\ast}_t, \mu^{j,\ast}_t, j \in \mfN, j \neq i\}$ given by \eqref{opt-dist}-\eqref{opt-mean}, and agent-$i$ is using an arbitrary distribution $\Phi^i_t \in \mc{U}^{i,\tex}$ with the mean value $\mu^i_t$, we have 
\begin{align}
&(i)~~ \mb{E} \Vert x_t^{i,\Phi^N} - x_t^{i,\Phi} \Vert^2 \leq C( o(\tfrac{1}{N}) + o (\delta_N^2)),\label{conv_order_state}\allowdisplaybreaks\\
&(ii)~~\Big\vert J_i^N(\Phi^{i},\Phi^{-i,\ast}) - J_i^{\infty}(\Phi^{i}) \Big\vert \leq C(o(\tfrac{1}{N}) + o(\delta_N)),\label{conv_order_cost}
\end{align}
where $C$ is a constant independent of $N$.$\hfill \square$
\end{thm}
\textit{Proof.} We prove each part separately. 
\textit{Part (i)}. From \eqref{dyn_modified_fin} and \eqref{dyn_modified_inf}, we have
\begin{align}
&d(x_t^{i,\Phi^N}-x_t^{i,\Phi}) = \Big(A_k (x_t^{i,\Phi^N}-x_t^{i,\Phi}) + (F_k \tfrac{1}{N}\sum_{j\neq i}x_t^{j,\Phi^N}-\bar{F}_k \bar{x}_t^{\Phi}) \nonumber\\&\hspace{0.1cm} + (H_k \tfrac{1}{N} \sum_{j\neq i}\mu^{j,N}_t -\bar{H}_k \bar{\mu}_t ) + F_k \tfrac{1}{N}x_t^{i,\Phi^N}+ H_k \tfrac{1}{N} \mu^{i}_t \Big)dt. 
\end{align}
The above ODE may be written as
\begin{align}\label{diff_state_fin_inf}
&x_t^{i,\Phi^N}-x_t^{i,\Phi} = \int_{0}^{t}e^{A_k (t-\tau)} (F_k \tfrac{1}{N}\sum_{j\neq i}x_\tau^{j,\Phi^N}-\bar{F}_k \bar{x}_\tau^{\Phi})d\tau \nonumber\\&\hspace{1.5cm} + \int_{0}^{t}e^{A_k (t-\tau)} (H_k \tfrac{1}{N} \sum_{j\neq i}\mu^{j,N}_\tau -\bar{H}_k \bar{\mu}_\tau)d\tau \nonumber\\&\hspace{1cm}+ \int_{0}^{t}e^{A_k (t-\tau)}(F_k \tfrac{1}{N}x_\tau^{i,\Phi^N}+ H_k \tfrac{1}{N} \mu^{i}_\tau)d\tau. 
\end{align}
Using $\Vert a+b+c \Vert^2 \leq 3 \Vert a\Vert^2 + 3 \Vert b\Vert^2 + 3 \Vert c\Vert^2$ for $a,b,c \in \mathbb{R}^n$, we get
\begin{multline}
\mb{E}\Vert x_t^{i,\Phi^N}\!\!\!-x_t^{i,\Phi} \Vert^2 \leq 3 \mb{E}\left\Vert \int_{0}^{t}e^{A_k (t-\tau)} (\tfrac{F_k}{N}\sum_{j\neq i}x_\tau^{j,\Phi^N}\!\!\!-\bar{F}_k \bar{x}_\tau^{\Phi})d\tau \right\Vert^2
\allowdisplaybreaks\\
+ 3 \mb{E} \left\Vert \int_{0}^{t}e^{A_k (t-\tau)} (H_k \tfrac{1}{N} \sum_{j\neq i}\mu^{j,N}_\tau -\bar{H}_k \bar{\mu}_\tau)d\tau \right\Vert^2
\allowdisplaybreaks\\
+ 3 \mb{E}\left\Vert \int_{0}^{t}e^{A_k (t-\tau)}(F_k \tfrac{1}{N}x_\tau^{i,\Phi^N}+ H_k \tfrac{1}{N} \mu^{i}_\tau)d\tau \right\Vert^2.
\end{multline}
Cauchy-Schwartz and Jensen's inequalities further imply
\begin{align}
&\mb{E}\Vert x_t^{i,\Phi^N}\!\!\!\!-x_t^{i,\Phi} \Vert^2 \leq 3 \!\!\! \int_{0}^{t}\Vert e^{A_k (t-\tau)}\Vert^2 \;\mb{E}\left\Vert \tfrac{F_k}{N}\sum_{j\neq i}x_\tau^{j,\Phi^N}\!\!\!\!-\bar{F}_k \bar{x}_\tau^{\Phi}\right\Vert ^2\!\!\!d\tau \nonumber\allowdisplaybreaks\\&\hspace{1.3cm} + 3 \int_{0}^{t} \Vert e^{A_k (t-\tau)} \Vert^2 \;\mb{E}\left\Vert H_k \tfrac{1}{N} \sum_{j\neq i}\mu^{j,N}_\tau -\bar{H}_k \bar{\mu}_\tau \right\Vert^2\!\!\! d\tau \nonumber\allowdisplaybreaks\\&\hspace{0.7cm}
+ 3 \int_{0}^{t} \Vert e^{A_k (t-\tau)} \Vert^2 \;\mb{E} \left\Vert F_k \tfrac{1}{N}x_\tau^{i,\Phi^N}+ H_k \tfrac{1}{N} \mu^{i}_\tau \right\Vert^2 d\tau.
\end{align}
Next, we compute the rate of convergence for each of the terms in the above inequality. For the first term we write
\begin{align}\label{ave_diff}
&\mb{E}\left\Vert F_k \tfrac{1}{N}\sum_{j\neq i}x_t^{j,\Phi^N}-\bar{F}_k \bar{x}_t^{\Phi} \right\Vert^2 = \mb{E}\left\Vert F_k \tfrac{1}{N}\sum_{j\neq i}(x_t^{j,\Phi^N}-\mb{E}x_t^{j,\Phi^N})
\right.\nonumber\allowdisplaybreaks\\&
\left.\hspace{1cm}+(F_k \tfrac{1}{N}\sum_{j\neq i}\mb{E}x_t^{j,\Phi^N}-F_k^N \bar{x}_t^{\Phi})+
(F_k^N -\bar{F}_k)\bar{x}_t^{\Phi}\right\Vert^2,
\end{align}
where $F^N := \pi^N \otimes F$, and 
\begin{align}\label{fin_ave}
F_k \tfrac{1}{N}\sum_{j\neq i}\mb{E}x_t^{j,\Phi^N} &= F_k \tfrac{1}{N} (N_1 \mb{E}x_t^{(.),1,\Phi^N}+ \dots + N_K \mb{E}x_t^{(.),K,\Phi^N}) \nonumber\\&-  F_k \tfrac{1}{N} \mb{E}x_t^{i,\Phi^N} = F^N_k \bar{x}_t^{\Phi} -  F_k \tfrac{1}{N} \mb{E}x_t^{i,\Phi^N}.
\end{align}
Substituting \eqref{fin_ave} in \eqref{ave_diff}, we have 
\begin{align}\label{ave_diff_2}
&\mb{E}\left\Vert F_k \tfrac{1}{N}\sum_{j\neq i}x_t^{j,\Phi^N}-\bar{F}_k \bar{x}_t^{\Phi} \right\Vert^2 \leq 3 \mb{E}\left\Vert F_k \tfrac{1}{\sqrt{N}}\sum_{j\neq i}\tfrac{(x_t^{j,\Phi^N}-\mb{E}x_t^{j,\Phi^N})}{\sqrt{N}} \right\Vert^2 \nonumber\\& \hspace{1.7cm}+3 \left\Vert F_k \tfrac{1}{N} \mb{E}x_t^{i,\Phi^N} \right\Vert^2+
3 \left\Vert (F_k^N -\bar{F}_k)\bar{x}_t^{\Phi}\right\Vert^2.
\end{align}
Since for $j \neq l$ the terms $(x_t^{j,\Phi^N}-\mb{E}x_t^{j,\Phi^N})$ and  $(x_t^{l,\Phi^N}-\mb{E}x_t^{l,\Phi^N})$ are independent, we can write 
\begin{align}
\mb{E}\Vert F_k \tfrac{1}{\sqrt{N}}\sum_{j\neq i}\tfrac{(x_t^{j,\Phi^N}-\mb{E}x_t^{j,\Phi^N})}{\sqrt{N}} \Vert^2 \leq C_1 \tfrac{1}{N} \mb{E} \sum_{j\neq i}\Vert \tfrac{(x_t^{j,\Phi^N}-\mb{E}x_t^{j,\Phi^N})}{\sqrt{N}} \Vert^2.
\end{align}
Hence 
\begin{equation}\label{ave_conv}
\mb{E}\Vert F_k \tfrac{1}{N}\sum_{j\neq i}x_t^{j,\Phi^N}-\bar{F}_k \bar{x}_t^{\Phi} \Vert^2 \leq C (O(\tfrac{1}{N}) + O(\delta_N^2)).
\end{equation}
Similarly we can show that 
\begin{equation}
\mb{E}\Vert H_k \tfrac{1}{N} \sum_{j\neq i}\mu^{j,N}_\tau -\bar{H}_k \bar{\mu}_t \Vert^2 \leq C(O(\tfrac{1}{N}) + O (\delta_N^2)),
\end{equation}
and hence \eqref{conv_order_state} follows.\\\\
\textit{Part (ii)}. For an arbitrary control distribution $\Phi^i$ and the optimal control distributions $\Phi^{j,*}, j\neq i$, we have 
\begin{multline}
J_i^N(\Phi^{i},\Phi^{-i,\ast}) - J_i^{\infty}(\Phi^{i}) = \mb{E} \bigg[\int_0^{\infty} e^{-\rho t}\bigg(\hf \Vert{x_t^{i,\Phi^N}}-y_t^{\Phi^N}\Vert^2_{Q_k} \allowdisplaybreaks\\ - \hf \Vert x_t^{i,\Phi}-\bar{y}^{\Phi}_t \Vert^2_{Q_k} + \eta_k^\trns (x_t^{i,\Phi^N} - {x_t^{i,\Phi}}) + \eta_k^\trns (\bar{y}_t^{\Phi}-y_t^{\Phi^N})\allowdisplaybreaks\\
+ \int_{\dsR^m}\Big(({x_t^{i,\Phi^N}} - {x_t^{i,\Phi}})^\trns + ({\bar{y}_t^{\Phi}} - y_t^{\Phi^N})^\intercal\Big)S_k u\Phi_t^{i}(u)du \bigg) dt \bigg],
\end{multline}
which reduces to
\begin{multline}
J_i^N(\Phi^{i},\Phi^{-i,\ast}) - J_i^{\infty}(\Phi^{i}) \allowdisplaybreaks\\
= \mb{E} \bigg[ \int_0^{\infty} e^{-\rho t}\bigg(\underbrace{\hf ({x_t^{i,\Phi^N}})^\trns Q_k {x_t^{i,\Phi^N}} - \hf ({x_t^{i,\Phi}})^\trns Q_k {x_t^{i,\Phi}}}_{D_1}\allowdisplaybreaks\\
+\underbrace{\hf (y_t^{\Phi^N})^\trns Q_k y_t^{{\Phi}^N} - \hf (\bar{y}_{t}^{\Phi})^\trns Q_k \bar{y}_t^{\Phi}}_{D_1^\prime}\allowdisplaybreaks\\
+ \underbrace{(x_t^{i,\Phi^N} - x_t^{i,\Phi})^\trns(\eta_k +Q_ky^{{\Phi}^N}_t)}_{D_2} - \underbrace{(\eta_k - Q_k x^{i,\Phi}_t)^\trns(y_t^{\Phi^N} - \bar{y}_t^{\Phi})}_{D_2^\prime}\allowdisplaybreaks\\+\underbrace{({x_t^{i,\Phi^N}} - {x_t^{i,\Phi}})^\trns S_k \mu_t^{i}}_{D_3}-\underbrace{(y_t^{\Phi^N} - \bar{y}_t^{\Phi})^\trns S_k \mu_t^{i}}_{D_3^\prime} \bigg) dt \bigg].
\end{multline}
The rate of convergence for $D_1$ is computed as
\begin{align}
&\mb{E} \left[D_1\right] = \hf \mb{E} \left[ \Vert {x_t^{i,\Phi^N}} - {x_t^{i,\Phi}} \Vert^2_{Q_k} \right] +  \mb{E} \left[ ({x_t^{i,\Phi}})^\trns Q_k({x_t^{i,\Phi^N}} - {x_t^{i,\Phi}}) \right] \nonumber\\ 
&\leq \mb{E} \left[ \Vert {x_t^{i,\Phi^N}} - {x_t^{i,\Phi}} \Vert^2_{Q_k} \right] + C\, \mb{E}\left[\Vert x_t^{i,\Phi}\Vert^2\right]^{\hf}\, \mb{E}\left[\Vert{x_t^{i,\Phi^N}} - {x_t^{i,\Phi}}\Vert^2\right]^{\hf} \nonumber\\
& = o(\delta_N) + o(\tfrac{1}{\sqrt{N}}).
\end{align}
Following the same approach and using \eqref{ave_conv}, we have 
\begin{equation}
\mb{E}\left[D^\prime_1 \right] = o(\delta_N) + o(\tfrac{1}{\sqrt{N}}).
\end{equation}
Moreover, from \eqref{conv_order_state}, we have
\begin{align}
\mb{E}\left[ D_2 \right] &= \mb{E} \left[(x_t^{i,\Phi^N} - x_t^{i,\Phi})^\trns(\eta_k +Q_k y^{{\Phi}^N}_t)\right] \nonumber\\
&\leq C\, \mb{E}\left[\Vert \eta_k +Q_k y^{{\Phi}^N}_t \Vert^2\right]^{\hf}\, \mb{E}\left[\Vert{x_t^{i,\Phi^N}} - {x_t^{i,\Phi}}\Vert^2\right]^{\hf} \nonumber\\
&= o(\delta_N) + o(\tfrac{1}{\sqrt{N}}).
\end{align}
Similarly we can show that 
\begin{equation*}
\mb{E}\left[ D_2^\prime \right] =\mb{E} \left[D_3\right] =\mb{E}\left[D_3^\prime\right]=o(\delta_N) + o(\tfrac{1}{\sqrt{N}}).\qquad\qquad\quad\hfill \square
\end{equation*}
Now we use the results of Theorem \ref{thm:fin_inf_diff} to show the $\epsilon$-Nash property for the solutions to the limiting exploratory LQG MFG.
\begin{thm}[$\epsilon$-Nash Property] \label{thm:eps_nash} The set of optimal distributions $\{\Phi_t^{i,\ast}, i \in \mfN\}$ with the corresponding mean values $\{\mu_t^{i,\ast}, i \in \mfN\}$ yields an $\epsilon-Nash$ equilibrium for the system \eqref{dyn_modified_fin}-\eqref{cost_exp_fin_pop}, i.e. 
\begin{equation*}
J_i^N(\Phi^{i,\ast},\Phi^{-i,\ast}) - \epsilon \leq \inf_{\Phi^{i}} J_i^N(\Phi^{i},\Phi^{-i,\ast}) \leq J_i^N(\Phi^{i,\ast},\Phi^{-i,\ast}).\,\,\,\hfill \square
\end{equation*}
\end{thm}
\textit{Proof.} From Theorem \ref{thm:fin_inf_diff}, we have
\begin{multline}\label{eps_ineq_arbitCntrl}
-o(\delta_N) - o(\tfrac{1}{\sqrt{N}})
 \leq J_i^N(\Phi^{i},\Phi^{-i,\ast}) - J_i^{\infty}(\Phi^{i})\allowdisplaybreaks\\ \leq o(\delta_N) + o(\tfrac{1}{\sqrt{N}}).
\end{multline}
Since \eqref{eps_ineq_arbitCntrl} holds for every $\Phi^{i}$, we can write
\begin{gather}
J_i^N(\Phi^{i,\ast},\Phi^{-i,\ast}) - o(\delta_N) - o(\tfrac{1}{\sqrt{N}}) \leq J_i^{\infty}(\Phi^{i,\ast})\label{eq1}
\allowdisplaybreaks\\
-o(\delta_N) - o(\tfrac{1}{\sqrt{N}})
 \leq \inf_{\Phi^{i}} J_i^N(\Phi^{i},\Phi^{-i,\ast}) - \inf_{\Phi^{i}} J_i^{\infty}(\Phi^{i}) \label{eq2}
\end{gather}
From \eqref{eq1} and \eqref{eq2}, we have 
\begin{align*}
 J_i^N(\Phi^{i,\ast},\Phi^{-i,\ast})& - 2o(\delta_N) - 2o(\tfrac{1}{\sqrt{N}}) \nonumber\allowdisplaybreaks\\&\hspace{0cm} \leq \inf_{\Phi^{i}} J_i^{\infty}(\Phi^{i}) \leq \inf_{\Phi^{i}} J_i^N(\Phi^{i},\Phi^{-i,\ast}).\quad\hfill \square
\end{align*}
\section{Classical versus Exploratory LQG MFGs}\label{sec:equivalence}
In this section we compare the solutions to the classical and exploratory LQG MFG systems in the infinite-population limit. First we establish the equivalence of the existence of solutions.
\begin{thm}(Existence Equivalence of Solutions to Classical \& Exploratory LQG MFGs) There exists a solution to \eqref{MinorDyn_inf}-\eqref{Minorcost_inf} if and only if there exists a solution to \eqref{dyn_modified_inf}-\eqref{cost_func_exploratory}. $\hfill \square$
\end{thm}
\textit{Proof.}
The solutions to the classical and exploratory LQG MFG systems are given, respectively, by \eqref{fixed_point_eqs_classic}-\eqref{opt-mean_classic} and \eqref{opt-dist}-\eqref{fixed_point_eqs}. Note that only the mean value $\mu^{i,\ast}_t$ (given by \eqref{opt-mean}) of \eqref{opt-dist} and its mean field limit $\bar{\mu}_t$ (given by \eqref{cntlMF_mean_exp}) appear in \eqref{dyn_modified_inf} and these processes, respectively, coincide with $u^{i,\ast}_t$ (given by \eqref{opt-mean_classic}) and $\bar{u}_t$ (given by \eqref{ubar_classic}) in \eqref{MinorDyn_inf}. Moreover by inspection we have
\begin{gather}
\mb{E}[x^{i,\ast}_t] = \mb{E}[x^{i,\Phi, \ast}_t],\allowdisplaybreaks\\
\mb{E}[(x^{i,\ast}_t)^\intercal x^{i,\ast}_t] = \mb{E}[(x^{i,\Phi,\ast}_t)^\intercal x^{i,\Phi,\ast}_t],
\end{gather}
for every $t\in \mfT$, where $x^{i,\ast}_t$ and $x^{i,\Phi, \ast}_t$ are the optimal trajectories induced by $u^{i,\ast}_t$ and $\Phi^{i,\ast}_t$, respectively. Therefore $\Phi^{i,\ast} \in \mc{U}^{i,\tex}$ if and only if $u^{i,\ast} \in \mc{U}^i$.
$\hfill \square$

Next, we establish that  the classical and exploratory systems share the same state and control mean field distributions. 
\begin{thm}[Equivalence of the mean field distributions] The state and the control mean field distributions for the classical and exploratory LQG MFGs governed by 
\eqref{MinorDyn_inf}-\eqref{Minorcost_inf} and \eqref{dyn_modified_inf}-\eqref{cost_func_exploratory}, respectively, are identical, given the same initial state distribution. 
\end{thm}
\textit{Proof.} In the limiting MFG, the state and the control mean field distributions of subpopulation $k$ coincide the distribution of the state and the control action of a generic agent in subpopulation $k$. Given the same initial state distribution in both the classical and the exploratory LQG MFGs, the processes $x^{i,\ast}_t$ and $x^{i,\Phi^\ast}_t$, $i \in \mc{I}_k$, are identical as they satisfy the same SDE given by \eqref{MinorDyn_inf} and \eqref{dyn_modified_inf}, where $u^{i,\ast}_t, \bar{u}_t,$ and $\mu^{i,\ast}_t, \bar{\mu}_t,$ are substituted, respectively. Hence the state mean field distribution is the same in both cases. We note that for the class of problems we are considering only the mean value $\bar{x}_t$ of the state mean field distribution enters in the dynamics and cost functionals. Moreover, for the classical case the control mean field distribution is given by the distribution of $u^{i,\ast}_t$ and for the exploratory case by the distribution of $\mu^{i,\ast}_t$, which again are the same processes. 
$\hfill \square$

Next, we investigate the cost of exploration in the infinite-population limit. The cost of exploration is defined as the difference in the original cost \eqref{Minorcost_inf} when the exploratory optimal control distribution \eqref{opt-dist} and the classical optimal control \eqref{opt-mean_classic} are applied. 
\begin{thm}[Cost of Exploration]\label{Cst_expl} Given the same initial state $x_0 \in \dsR^n$ for a generic agent in sub-population $k \in \mathfrak{K}$ in both the classical and the exploratory cases, the cost of exploration $COE^{\infty,k}$ for the agent in the infinite-population limit is given by 
\begin{equation*}
\qquad\qquad\qquad\quad COE^{\infty,k}(x_0) = \frac{\lambda_k}{2\rho}. \qquad\qquad\qquad\qquad \hfill \square
\end{equation*}
\end{thm}
\textit{Proof.}
The cost of exploration is defined as in 
\begin{align}
&COE^{\infty,k}(x_0) = \bigg(J_i^\infty(\Phi^{i,\ast}) - \lambda_k\mb{E} \bigg[ \int_0^{\infty} e^{-\rho t}\nonumber\\ &\hspace{0.7cm}\times\Big(\int_{\dsR^m} \Phi_t^{i,\ast}(u)\ln\Phi_t^{i,\ast}(u) du\Big) dt\bigg| x^i_0 = x_0\bigg]\bigg)- J^{\infty}_i(u^{i,\ast}),
\end{align}
where the first term on the right hand side represents the original cost subject to the optimal control  distribution obtained in the exploratory version of the problem (it can be seen as the average cost when the applied actions are sampled from $\Phi^{i,*}_t$ (see \cite{XYZRL2020})).

From \eqref{dyn_modified_inf}, only the mean value $\mu^{i,\ast}_t$ of the optimal distribution $\Phi^{i,\ast}_t$ and its mean field limit $\bar{\mu}_t$ appear in the exploratory dynamics which happen to respectively coincide with $u^{i,\ast}_t$  and $\bar{u}_t$ in the classical dynamics \eqref{MinorDyn_inf}. Hence, we have $x_t^{i,\ast} = x_t^{i,\Phi^\ast}$, and
\begin{align}
    \mb{E} \bigg[\int_0^{\infty}& e^{-\rho t}\int_{\dsR^m}\Big(({x_t^{i,\Phi^\ast}})^\trns S_k + n_k^\trns \Big)u\Phi_t^{i,\ast}(u)du \bigg) dt \bigg] \nonumber\allowdisplaybreaks\\&
    = \mb{E} \bigg[\int_0^{\infty} e^{-\rho t}\Big(({x_t^{i,\ast}})^\trns S_k + n_k^\trns \Big)u^{i,\ast}_t dt \bigg].
\end{align}
Therefore, we have
\begin{multline}
COE^{\infty,k}(x_0)= \mb{E} \bigg[\int_0^{\infty} e^{-\rho t}\bigg(\int_{\dsR^m}\hf u^\trns R_k u \Phi_t^{i,\ast}(u)du \nonumber\allowdisplaybreaks\\\hspace{0.0cm}- \hf{(u_t^{i,\ast})}^\trns R_k {u_t^{i,\ast}}\bigg) dt \bigg \vert x^i_0 = x_0  \bigg]\nonumber\allowdisplaybreaks\\
= \mb{E} \bigg[\int_0^{\infty} e^{-\rho t}\!\!\!\int_{\dsR^m}\!\!\hf (u-u^{i,\ast}_t)^\trns R_k (u-u^{i,\ast}_t) \Phi_t^{i,\ast}(u)du  dt \bigg\vert x^i_0 = x_0 \bigg]\nonumber\allowdisplaybreaks\\
= \frac{\lambda_k}{2\rho}\hspace{0.0cm} \hspace{0.0cm}\hfill \square\nonumber
\end{multline}
We note that the mean value $\mu^{i,\ast}_t$ of the optimal Gaussian distribution $\Phi^{i, \ast}_t$ coincides with the optimal action $u^{i,\ast}_t$ for the classical LQG MFG problems where no exploration is involved. The mean value $\mu^{i,\ast}_t$ is a function of the agent-$i$'s state, and is not impacted by the exploration weight $\lambda_k$. The variance $\lambda_k R^{-1}_k$ of $\Phi^{i,\ast}_t$, in contrast, is directly impacted by $\lambda_k$. In fact when there is no exploration, the distribution $\Phi^{i,\ast}_t$ reduces to the Dirac delta probability mass. In the following theorem we formally show that the classical LQG MFG can be formulated as a special case of exploratory LQG MFG when the exploration weight $\lambda_k \rightarrow 0$, i.e. when there is no exploration. 
\begin{thm}\label{thm:convergence} Given an initial state $x_0 \in \dsR^n$ for agent-$i$, we have 
\begin{equation}
    (i)~~ \lim_{\lambda \rightarrow 0} \Phi_t^{i,\ast}(u) = \delta({u-u^{i,\ast}_t}),
\end{equation}
where $\delta$ denotes the Dirac delta function, and
\begin{gather}
  (ii)~~  \lim_{\lambda \rightarrow 0}|V_i^{\infty}(x_0)-V_i^{\tex,\infty}(x_0)|=0, 
\end{gather}
where $V_i^{\infty}(x_0)$ and $V^{\tex,\infty}_i(x_0)$, respectively, denote the classical and exploratory value functions for agent-$i$ in the infinite-population case.
\end{thm}
\textit{Proof.} We have 
\begin{equation}
 \Phi^{i, \ast}_t = \mathcal{N}(\mu^{i,\ast}_t, \lambda_k R_k^{-1}),   
\end{equation}
where $\mu^{i,\ast}_t = u^{i,\ast}_t$ given by \eqref{opt-mean} is independent of $\lambda_k$. Hence as $\lambda_k \rightarrow 0$, $\Phi^{i, \ast}_t$ converges weekly to $\delta({u-u^{i,\ast}_t})$.

Moreover, substituting $\Phi^{i,\ast}_t$ in \eqref{cost_func_exploratory} results in 
\begin{multline}
    V_i^{\tex, \infty}(x_0) = J_i^{\infty}(\Phi^{i,\ast})\allowdisplaybreaks\\
     = \mb{E} \bigg[ \int_0^{\infty} e^{-\rho t}\bigg( \hf\Vert x_t^{i,\Phi}-\bar{y}^{\Phi}_t\Vert_{Q_k}^2     +\eta_k^\trns (x_t^{i, \Phi}-\bar{y}^{\Phi}_t)\allowdisplaybreaks \\
      + (x_t^{i,\Phi}-\bar{y}^{\Phi}_t)^\trns S_k u^{i,\ast}_t  + n_k^\trns u^{i,\ast}_t \allowdisplaybreaks\\
    + \int_{\dsR^m}\Big(\hf u^\trns R_k u+\lambda_k \ln \Phi_t^{i,\ast}(u) \Big)\Phi_t^{i,\ast}(u)du \bigg) dt \bigg\vert x^i_0 = x_0\bigg].
\end{multline}
Subsequently
\begin{multline}
V_i^{\infty}(x_0)-V_i^{\tex,\infty}(x_0) \\= \mb{E} \bigg[\int_0^{\infty} e^{-\rho t}\bigg(\hf{(u_t^{i,\ast})}^\trns R_k {u_t^{i,\ast}}-\int_{\dsR^m}\hf u^\trns R_k u \Phi_t^{i,\ast}(u)du \\- \lambda_k \int_{\dsR^m}\Phi_t^{i,\ast}(u)\ln \Phi_t^{i,\ast}(u)du\bigg) dt \bigg\vert x^i_0=x_0\bigg]\nonumber\allowdisplaybreaks\\
= \frac{\lambda_k}{2\rho}\big(\ln(2\pi \lambda_k R^{-1}_k)-1\big).
\end{multline}
Finally we have $\lim_{\lambda_k \rightarrow 0}| \frac{\lambda_k}{2\rho}\big(\ln(2\pi \lambda_k R^{-1}_k)-1\big)| = 0.$
$\hfill \square$

\section{Potential Applications in Learning MFGs}\label{sec:apps}

The randomisation around the optimal control action generated through sampling according to the Gaussian distribution \eqref{opt-dist} could potentially be used for obtaining consistent estimation algorithms in the context of adaptive control. Moreover, the Gaussian distribution could be used as a policy approximation in model-based reinforcement learning (RL) in LQG MFG systems for learning the dynamical parameters or the value (reward) function directly (for an example of the latter in a single-agent mean-variance problem see \cite{XYZ-MV-2020}). Furthermore, the exploratory LQG MFGs can be used to approximate nonlinear MFGs which are not explicitly solvable in general. This allows to explore the states of nonlinear system which are not captured by the classical LQG MFG approximations and is expected to improve the performance of learning algorithms for such systems. Below we suggest a preliminary learning algorithm for the setup introduced in the paper in the context of optimal executions in electronic markets. 




In electronic markets traders typically know their cost functional, but are unsure of the market model and their impact on the price. Consider the case where there are $N$ high-frequency traders in the market. For trader $i \in \mfN$, the trading dynamics is given by  
\begin{subequations}
\label{tradingDynamics}
\begin{align}
dF_t &= \lambda\,\, \frac{1}{N} \sum_{j=1}^{N} \nu^j_t dt + \sigma\, dw_t,
& 
dq_t^i &= \, \nu_t^idt,
\\
S_t^i &= F_t+a\,\int_0^t\nu_u^i\,du,
&
dZ_t^i &= -S^i_t \, dq_t^i,
\end{align}
\end{subequations}
where $q_t^i, \nu^i_t, S_t^i, Z^i_t$ denote, respectively, the inventory, trading rate, execution price, and cash process for agent-$i$ at time $t \in \mfT$. Moreover, $F_t$ denotes the (fundamental) midprice and $w_t$ denotes a Wiener process modeling the impact of noise (uninformed) traders on price. The unknown parameters to be learned in this model are the volatility $\sigma>0$, the permanent price impact strength $\lambda\ge0$, and the temporary price impact strength $a\ge0$.

The cost functional for agent-$i$ is given by
\begin{equation}\label{tradingCost}
    J(\nu^i, \nu^{-i}) = \mathbb{E}\biggl[\tfrac{\phi}{2}  \int_0^T (q_u^i)^2 du - Z_T^i -q_T^i \big(F_T - \psi q_T^i\big)\biggr],
\end{equation}
where the running and terminal costs correspond, respectively, to an urgency penalty and terminal book value. The parameters $\phi\ge0$ and $\psi\ge0$ denote the urgency parameter and the terminal execution penalty parameter, respectively.  The optimization problem \eqref{tradingDynamics}-\eqref{tradingCost} can be formulated as an LQG MFG problem \cite{FirooziISDG2017,casgrain2018mean}, and is indeed a special case of the general system described by \eqref{MinorDyn}-\eqref{minorCost}. Hence the optimal control distribution is given by \eqref{opt-dist}.
Subsequently an algorithm for learning the dynamical parameters can be designed. The algorithm consists of four phases detailed below, mainly, initialization, learning, planning, and acting.
\begin{table}[h]
\begin{tabular}{l}
\hline
\textbf{RL Algorithm for Agent-$\boldsymbol{i}$}\\ \hline
\begin{minipage}[t]{0.4\textwidth} 
\begin{itemize}
    \item \textit{Initialization}: Run the base policy $\Phi_0^{i,\ast}(u)$. This policy could be obtained using the optimal control distribution \eqref{opt-dist} with an arbitrary/initial choice of model parameters.   
\item \textit{Model Learning}: Estimate the dynamical parameters using the new set of data collected (e.g. maximum likelihood could be used). This step involves estimating the parameters in the mean field equation. 
\item\textit{Planning}: Update the policy $\Phi_t^{i,\ast}(u)$ using the updated value of the parameters or the state in \eqref{opt-dist}.
\item \textit{Acting}:
\begin{itemize}
    \item[(i)] Sample an action from the distribution $\Phi_t^{i,\ast}(u)$.
    \item[(ii)] Execute the action. This involves applying the action to equation (3), where the current values of the parameters and the mean fields are used. 
    \item[(iii)] Add the resulting data $\{(x_t^i, u_t^i)\}$ to the data set $\mathcal{D}$.
\end{itemize}
\item Repeat \textit{Planning} and \textit{Acting} steps $n$ times. 
\item Go back to \textit{Model Learning}
\end{itemize} 
\end{minipage}\\ \hline
\end{tabular}
\end{table}

\section{Concluding Remarks}\label{sec:conclusion}

Theorem \ref{thm:convergence} shows that in the limit as $\lambda \rightarrow 0$, the solutions of exploratory LQG MFG converges to that of classical LQG MFGs (i.e., there is no exploration). Moreover, the mean value of the optimal exploratory distribution coincides with the classical optimal control  (analogous to the single-agent case in \cite{XYZRL2020}). This suggests that exploitation and exploration are, respectively, reflected by the mean value and the variance of the optimal distribution. Therefore, for exploration we do not need to consider all possible control distributions. The best distribution to balance exploration versus exploitation for the LQG MFG setup considered in this work is Gaussian according to Theorem \ref{thm:opt_dist}. This could be considered as a refinement of the $\epsilon$-greedy policy which is a common approach to balance exploitation versus exploration, especially in  multi-armed bandit problems, where a Bernoulli distribution is used such that the current best arm is selected with a probability $1-\epsilon$ (exploit), and otherwise a random arm is selected with a probability $\epsilon$ (explore).

\bibliographystyle{plain}        
\bibliography{references}           
                                

\appendix  
\section{Notations}
The important notations that are used in the paper are summarized below for the conveniece of the reader.   
\begin{itemize}
    \item $\mc{I}_k = \lbrace i : \theta_i = \theta^{(k)},~ i \in \mfN \rbrace$,
    \item $\Vert a\Vert^2_{R}=a^\trns R a$,
    \item $u^{-i} = (u^{1}, \dots, u^{i-1}, u^{i+1}, \dots u^{N})$,
    \item  $\pi_k^N = \tfrac{N_k}{N}, ~ k \in \mfK$,
    \item $\pi^{N} = (\pi_{1}^{N},...,\pi_{K}^N)$,
    \item $\pi_k = \lim_{N\rightarrow \infty}\tfrac{N_k}{N}$,
    \item $\pi = (\pi_{1},...,\pi_{K})$,
   \item $\bar{F}_k:= F_k \otimes \begin{bmatrix}\pi_1,...,\pi_K\end{bmatrix}$, 
   \item $\bar{H}_k:= H_k \otimes \begin{bmatrix}\pi_1,...,\pi_K\end{bmatrix}$,
   \item $\bar{\psi}_k:= \psi_k \otimes \begin{bmatrix}\pi_1,...,\pi_K\end{bmatrix}$.
\end{itemize}
\subsection{Finite-Population Processes}
\subsubsection{Classical Case}
\begin{itemize}
    \item $x^i_t$: state of agent-$i$ at time $t$,
    \item $u^i_t$: control of agent-$i$ at time $t$,
    \item $J^{N}_i(u^i, u^{-i})$: cost functional for agent $i$, 
    \item $x^{(N)}:= \frac{1}{N} \sum_{i \in \mfN} x^i_t$,
    \item $u^{(N)}_t := \frac{1}{N} \sum_{i \in \mfN} u^i_t$,
\end{itemize}

\subsubsection{Exploratory Case}
\begin{itemize}
    \item $\Phi^{i,N}_t$: control distribution of agent-$i$ at time $t$,
    \item $J_i(\Phi^{i,N},\Phi^{-i,N})$: cost functional for agent $i$,
    \item $x_t^{i,\Phi^N}$:  state of agent-$i$ with actions sampled wrt $\Phi^{i,N}_t$,
    \item $\mu_t^{i,N} = \int_{\mathds{R}^m} u \Phi_t^{i,N}(u)du$,
    \item $x_t^{(N),\Phi^N} = \frac{1}{N} \sum_{i\in \mfN} x^{i,\Phi^N}_t$,
    \item $\mu_t^{(N)} = \frac{1}{N} \sum_{i \in \mfN} \mu^{i,N}_t$,
    \item $y_t^{\Phi^N}=\psi_k x_t^{(N),\Phi^N}$.
\end{itemize}
\subsection{Infinite-Population Processes}
\subsubsection{Classical Case}
\begin{itemize}
     \item $x^i_t$: state of agent-$i$ at time $t$,
     \item $u^i_t$: control of agent-$i$ at time $t$,
     \item $J^{\infty}_i(u^i)$: cost functional for agent $i$, 
      \item $\bar{x}_t^k := \lim_{N_k\rightarrow \infty} \frac{1}{N_k}\sum_{i\in \mc{I}_k} x^{i,k}_t$,
    \item $\bar{u}_t^{k} := \lim_{N_k\rightarrow \infty} \frac{1}{N_k}\sum_{i\in \mc{I}_k} u^{i,k}_t$,
    \item $(\bar{x}_t)^\trns = \begin{bmatrix} (\bar{x}_t^{1})^{\trns}, \dots,  
     (\bar{x}_t^{K})^{\trns}\end{bmatrix}$,
   \item $\bar{u}_t^\trns = \begin{bmatrix} (\bar{u}_t^1)^{\trns}, \dots, (\bar{u}_t^K)^\trns\end{bmatrix}$,
   \item $\bar{y}_t = \bar{\psi_k} \bar{x}_t$.
\end{itemize}

\subsubsection{Exploratory Case}
\begin{itemize}
    \item $\Phi^{i}_t$: control distribution of agent-$i$ at time $t$,
    \item $J_i^\infty(\Phi^i)$: cost functional for agent $i$,
    \item $\Phi_t^{i,\ast}(u)$: optimal control distribution of agent-$i$ at time $t$, 
    \item $x_t^{i,\Phi}$: state of agent-$i$ with actions sampled wrt $\Phi^{i}_t$,
    \item $\mu^i_t = \int_{\mathds{R}^m} u \Phi_t^{i}(u)du$,
    \item $\mu_t^{i,\ast}=\int_{\mathds{R}^m} u \Phi_t^{i,\ast}(u)du$,
    \item $\bar{x}_t^{k,\Phi} := \lim_{N_k \rightarrow \infty} \frac{1}{N_k} \sum_{i\in \mc{I}_k} x^{i,\Phi}_t$,
    \item $(\bar{x}_t^{\Phi})^\trns = \begin{bmatrix} (\bar{x}_t^{1,\Phi})^{\trns}, \dots, (\bar{x}_t^{K,\Phi})^{\trns}\end{bmatrix}$, 
    \item $\bar{\mu}_t^\trns = \begin{bmatrix} (\bar{\mu}_t^1)^{\trns}, \dots, (\bar{\mu}_t^K)^\trns\end{bmatrix}$, 
    \item $\bar{\mu}_t^k :=\lim_{N_k \rightarrow \infty} \frac{1}{N_k} \sum_{i \in \mc{I}_k} \mu^i_t$,
    \item $\bar{y}_t^{\Phi} = \bar{\psi}_k \bar{x}_t^{\Phi}$.
\end{itemize}



\end{document}